\documentclass[a4paper,11pt]{article}
\usepackage[utf8]{inputenc}
\usepackage{enumerate}
\usepackage{lmodern}
\usepackage[T1]{fontenc}
\usepackage{verbatim}
\usepackage{textcomp}
\usepackage[english]{babel}
\usepackage[a4paper,vmargin={3.5cm,3.5cm},hmargin={2.5cm,2.5cm}]{geometry}
\usepackage[font=sf, labelfont={sf,bf}, margin=1cm]{caption}
\usepackage[pdftex]{hyperref}

\usepackage[pdftex]{color,graphicx}

\usepackage{amsmath,amsfonts,amssymb,amsthm,mathrsfs}

\newcommand{\op}[1]{\operatorname{#1 }}

\newtheorem{theorem}{Theorem}

\newtheorem{lemma}[theorem]{Lemma}

\newtheorem{corollary}[theorem]{Corollary}
\newtheorem{proposition}[theorem]{Proposition}

\def\llbracket{[\hspace{-.10em} [ }
\def\rrbracket{ ] \hspace{-.10em}]}
\def\pp{{\mathcal P}}
\def\t{{\mathcal T}}
\def\ve{\varepsilon}
\def\wt{\widetilde}
\def\bm{{\bf m}}
\def\be{{\bf e}}
\def\bp{{\bf p}}
\def\la{\longrightarrow}
\def\R{{\mathbb R}}
\def\Q{{\mathbb Q}}
\def\SS{{\mathbb S}}
\def\ov{\overline}
\def\build#1_#2^#3{\mathrel{
\mathop{\kern 0pt#1}\limits_{#2}^{#3}}}
\def\rem{\noindent{\bf Remark. }}

\title{The Brownian plane}
\author{Nicolas Curien and Jean-Fran\c cois Le Gall}

\date{\small \it Ecole normale sup\'erieure and Universit\'e Paris-Sud}
\begin{document}
\maketitle

\begin{abstract}
We introduce and study the random non-compact metric space called the Brownian plane,
which is obtained as the scaling limit of the uniform infinite planar quadrangulation.
Alternatively, the Brownian plane is identified as the Gromov-Hausdorff tangent cone in distribution of the Brownian map
at its root vertex, and it also arises as the scaling limit of uniformly distributed (finite) planar quadrangulations with $n$ faces when the scaling factor tends to $0$ less fast than $n^{-1/4}$. We discuss various 
properties of the Brownian plane. In particular, we prove that the Brownian plane is homeomorphic to the plane, and
we get detailed information about geodesic rays to infinity. 
\end{abstract}

\section{Introduction}

A planar map is a finite connected (multi)graph drawn on the two-dimensional sphere and viewed up
to orientation-preserving homeomorphisms of the sphere. The faces of a planar
map are the connected components of the complement of edges, and the degree 
of a face counts the number of edges in its boundary, with the special 
convention that if both sides of an edge are incident to the same face then this edge is counted
twice in the degree of this face. Important special cases of planar maps are triangulations, where all faces have
degree three, and quadrangulations, where all faces have degree four. It is usual 
to deal with rooted planar maps, meaning that there is a distinguished edge, which
is also oriented and whose origin is called the root vertex.

Much recent work has been devoted to understanding the asymptotic properties 
of large planar maps chosen uniformly at random in a particular class, e.g. the class
of all triangulations, or of all quadrangulations, with a fixed number $n$ of faces tending to infinity. There are 
two basic kinds of limit theorems giving information about large random planar maps.

First, local limit theorems consider for every fixed integer $r\geq 1$ the combinatorial ball of radius 
$r$ in the planar map (this is the new planar map obtained by keeping only those faces whose boundary
contains at least one vertex whose graph distance from the root vertex is smaller than $r$), and show that
this ball converges in distribution as the size of the map tends to infinity towards
the corresponding ball in a random infinite planar lattice. Such local limits were studied 
first in the case of triangulations by Angel and Schramm \cite{Ang,AS} and the limiting object
is called the uniform infinite planar triangulation (UIPT). The analogous result for
quadrangulations was obtained later by Krikun \cite{Kri} (see also Chassaing and 
Durhuus \cite{CD} and M\'enard \cite{Men}), leading to the uniform infinite planar
quadrangulation (UIPQ). 

Secondly, scaling limits consist in looking at the vertex set of a planar map with $n$ faces as a metric space
for the graph distance rescaled by the factor $n^{-1/4}$, 
and studying the convergence of this metric space when $n$ tends 
to infinity, in the sense of the Gromov-Hausdorff distance familiar to geometers. The factor $n^{-1/4}$ is chosen 
so that the diameter of the rescaled planar map remains bounded in probability: It was first
noticed by Chassaing and Schaeffer \cite{CS} that the diameter of a random quadrangulation with $n$ faces
is of order $n^{1/4}$, and a similar result holds for much more general random planar maps,
including triangulations. The existence of a scaling limit for (uniformly distributed)
random quadrangulations was obtained recently in the papers \cite{LGU, Mi-quad}, leading to
a limiting random compact metric space called the Brownian map. In \cite{LGU}, it is also
proved that the Brownian map is the universal scaling limit of more general random planar maps
including triangulations. 

Our main goal in the present work is to provide a connection between the preceding 
limit theorems, by introducing a random (non-compact) metric space called the 
Brownian plane, which can be viewed either as the scaling limit of the UIPQ or 
as the Gromov-Hausdorff tangent cone in distribution of the Brownian map at its root. The Brownian plane can also 
be obtained as the limit of rescaled random quadrangulations with $n$ faces if the graph distance is multiplied by a
factor $\varepsilon_n$ such that $\varepsilon_n\to 0$ and $\varepsilon_nn^{1/4}\to\infty$ as $n\to\infty$.

Let us give a precise definition of the Brownian plane before stating our main results. We consider two independent
three-dimensional Bessel processes $R$ and $R'$ started from $0$ (see e.g. \cite{RY} for basic facts about Bessel processes). We then 
define a process $X=(X_t)_{t\in \R}$ indexed by the real line, by setting
$$X_t=\left\{\begin{array}{ll} R_t \qquad&\hbox{if }t\geq 0,\\
\noalign{\smallskip}
R'_{-t}\qquad&\hbox{if }t\leq 0.
\end{array}
\right.
$$
Then, for every $s,t\in\R$, we set
$$\ov{st} =\left\{
\begin{array}{ll}
[s\wedge t,s\vee t]\qquad&\hbox{if } st\geq 0,\\
\noalign{\smallskip}
(-\infty,s\wedge t] \cup [s\vee t,\infty)\qquad&\hbox{if } st< 0,
\end{array}
\right.$$
and
$$m_X(s,t)= \inf_{r\in \ov{st}} X_r.$$
We define a random pseudo-distance on $\R$ by 
$$d_X(s,t)= X_s + X_t -2\,m_X(s,t)$$
and put $s\sim_X t$ if $d_X(s,t)=0$. The quotient space $\t_\infty=\R/\!\sim_X$
equipped with $d_X$ is a (non-compact) random real tree, which is
sometimes called the infinite Brownian tree. This tree corresponds to Process 2 in
Aldous \cite{Al1}.  It can be verified 
that $\t_\infty$ is the tangent cone in distribution of Aldous' CRT 
at its root vertex, in a sense that will be explained below (see Theorem 11 in \cite{Al1} for a closely related result). 
We write $p_\infty:\R\la \t_\infty$
for the canonical projection and 
set $\rho_\infty= p_\infty(0)$, which plays the role of the root of $\t_\infty$.

We next consider Brownian motion indexed by $\t_\infty$. Formally, we consider
a real-valued process $(Z_t)_{t\in \R}$ such that,
conditionally given the process $X$, $Z$ is a centered Gaussian process  with
covariance
$$E[Z_sZ_t\mid X] =m_X(s,t)$$
so that we have $Z_0=0$ and $E[(Z_s-Z_t)^2\mid X] = d_X(s,t)$. It is not hard to verify that
the process $Z$ has a modification with continuous paths, and we consider this
modification from now on. Then a.s. we have $Z_s=Z_t$ for every $s,t\in \R$
such that $d_X(s,t)=0$, and thus we may (and sometimes will) view $Z$
as indexed by $\t_\infty$.

For every $s,t\in \R$, we set
$$D^\circ_\infty (s,t) = Z_s + Z_t - 2\,\min_{r\in[s\wedge t,s\vee t]} Z_r.$$
We extend the definition of $D^\circ_\infty$ to $\t_\infty\times \t_\infty$
by setting for $a,b\in\t_\infty$,
$$D^\circ_\infty (a,b)= \min\{ D^\circ_\infty (s,t): s,t\in\R,\; p_\infty(s)=a, p_\infty(t)=b\}.$$
Finally, we set, for every $a,b\in\t_\infty$,
$$D_\infty(a,b) = \inf_{a_0=a,a_1,\ldots,a_p=b} \sum_{i=1}^p D^\circ_\infty(a_{i-1},a_i)$$
where the infimum is over all choices of the integer $p\geq 1$ and of the
finite sequence $a_0,a_1,\ldots,a_p$ in $\t_\infty$ such that $a_0=a$ and
$a_p=b$. It is not hard to verify that $D_\infty$ is a pseudo-distance on $\t_\infty$.
We put $a\approx b$ if $D_\infty (a,b)=0$ (as will be show in Proposition \ref{identrel}, this is equivalent 
to the property $D^\circ_\infty (a,b)=0$). The Brownian plane is the quotient space $\pp=\t_\infty /\approx$, which is
equipped with the metric induced by $D_\infty$ and with the distinguished 
point which is the equivalence class of $\rho_\infty$. We simply write 
$\rho_\infty$ for this equivalence class, and use the notation $\boldsymbol{ \mathcal{P}}=(\pp, D_\infty,\rho_\infty)$ for the Brownian plane
viewed as a pointed metric space (recall that a metric space $(E,d)$ is said to be pointed if there
is a distinguished point $\alpha\in E$).

In order to state our first result,
we introduce the notation $(\mathrm{m}_\infty,D^*)$ for the Brownian map, as defined 
in the introduction of \cite{LGU} or in \cite{Mi-quad} for instance. Recall that $\mathrm{m}_\infty$ is obtained 
as a quotient set of Aldous' CRT \cite{Al1,Al3}, and that in this construction, the Brownian 
map comes with a distinguished point, which is the equivalence class of 
the root $\rho$ of the CRT. We will abuse 
notation and also write $\rho$ for the equivalence class of $\rho$ in $\mathrm{m}_\infty$.
From our perspective, it will be important to view the Brownian map as
a triplet $\mathbf{m}_\infty:=(\mathrm{m}_\infty, D^*,\rho)$, which is
a (random) pointed compact metric space. Note that, in a sense 
that can be made precise, $\rho$ is uniformly distributed over $\mathrm{m}_\infty$.

For every $\ve>0$, let $B_\ve(\boldsymbol{ \mathcal{P}})$ be the closed ball of radius $\ve$ centered at $\rho_\infty$ in $\mathcal{P}$,
and let $B_\ve(\mathbf{m}_\infty)$ be the closed ball of radius $\ve$ centered at $\rho$ in $\mathrm{m}_\infty$. Each of these balls 
is pointed at its center and thus viewed as a pointed metric space.

\begin{theorem}
\label{tangentcone}
For every $\delta >0$, we can find $\ve>0$ and construct on the same probability space
copies of the Brownian plane $\boldsymbol{ \mathcal{P}}$ and of the Brownian map
$\mathbf{m}_\infty$, in such a way that the balls $B_\ve(\boldsymbol{ \mathcal{P}})$
and  $B_\ve(\mathbf{m}_\infty)$ are isometric with probability at least $1-\delta$. 
Furthermore, we have
\begin{equation}
\label{D0}
(\mathrm{m}_\infty, \lambda D^*, \rho)\build{\la}_{\lambda\to\infty}^{(d)} (\pp, D_\infty,\rho_\infty)
\end{equation}
in distribution for the local Gromov-Hausdorff topology.
\end{theorem}

Let us briefly discuss the local Gromov-Hausdorff topology (for more details see 
Chapter 8.1 in \cite{BBI} and Section 2 below). First recall that a metric space $(E,d)$
is called a length space if, for every $a,b\in E$, the distance $d(a,b)$ coincides with the infimum
of the lengths of continuous curves connecting $a$ to $b$. We also say that $(E,d)$ is boundedly compact
if all closed balls are compact.
Then a sequence $(E_n,d_n,\alpha_n)$
of pointed boundedly compact length spaces is said to converge to $(E,d,\alpha)$
in the local Gromov-Hausdorff topology
if, for every $r>0$, the closed ball of radius $r$ centered at $\alpha_n$
in $E_n$ converges to the closed ball of radius $r$ centered at $\alpha$ in $E$, for the
usual Gromov-Hausdorff distance between pointed compact metric spaces. The space of all
pointed boundedly compact length spaces (modulo isometries) can be equipped with a
metric which is compatible with the preceding notion of convergence, and is then
separable and complete for this metric. The convergence in Theorem \ref{tangentcone} is just 
the usual convergence in distribution for random variables with values in
a Polish space. 

If $(E,d)$ is a boundedly compact length space and $\alpha\in E$, the Gromov-Hausdorff tangent cone
of $(E,d)$ at $\alpha$, if it exists, is the limit in the local Gromov-Hausdorff topology 
of the rescaled spaces $(E,\lambda d,\alpha)$ when $\lambda$ tends 
to infinity (see Section 8.2 in \cite{BBI}). The convergence (\ref{D0}) can thus be interpreted 
by saying that the Brownian plane is the Gromov-Hausdorff tangent cone in distribution of the
Brownian map at its root.

\medskip

\rem It follows from the convergence (\ref{D0}) that the Brownian plane is scale invariant, meaning that,
for every $\lambda >0$, $(\pp, \lambda D_\infty,\rho_\infty)$ has the same distribution as 
$(\pp, D_\infty,\rho_\infty)$. This can also be verified directly from the definition, using the similar
property for the infinite Brownian tree. 

\medskip
In order to state our second theorem, let us write $Q_\infty$ for the  uniform infinite planar quadrangulation
and $V(Q_\infty)$ for the vertex set of $Q_\infty$. The root vertex of $Q_\infty$
is denoted by $\rho_{(\infty)}$. Similarly, for every integer $n\geq 1$, $Q_n$ stands for a uniformly distributed 
rooted quadrangulation
with $n$ faces, $V(Q_n)$ is the vertex set of $Q_n$ and $\rho_{(n)}$ is the root vertex of 
$Q_n$. Finally $\op{d}_{\rm gr}$ denotes the graph distance on either $V(Q_n)$
or $V(Q_\infty)$.

\noindent \begin{theorem} \label{scaling-limit} We have 
\begin{eqnarray}
(V(Q_\infty), \lambda \op{d}_{\rm gr}, \rho_{(\infty)})& \xrightarrow[\lambda\to 0]{(d)} & (\pp, D_\infty,\rho_\infty), \label{D2}  \end{eqnarray}
in distribution for the local Gromov-Hausdorff topology. Furthermore, 
let $(k_{n})_{n\geq0}$ be a sequence of non-negative real numbers such that $k_{n} \to \infty$ and  $k_{n} = o(n^{1/4})$ as $n$ tends to infinity. Then, \begin{eqnarray} (V(Q_n), k_{n}^{-1}\op{d}_{\rm gr}, \rho_{(n)}) & \xrightarrow[n\to\infty]{(d)} &  (\pp, D_\infty,\rho_\infty)\label{D1}  \end{eqnarray}
in distribution for the local Gromov-Hausdorff topology.
\end{theorem}

The reader may have noticed that neither $(V(Q_\infty), \op{d}_{\rm gr})$ nor $(V(Q_n), \op{d}_{\rm gr})$ is a length space, so that
the discussion following Theorem \ref{tangentcone} does not seem to apply to the convergences (\ref{D2}) and (\ref{D1}). However it is
very easy to approximate $(V(Q_\infty), \op{d}_{\rm gr})$ (resp. $(V(Q_n), \op{d}_{\rm gr})$) by a pointed boundedly compact length space,
in such a way that the balls centered at the distinguished point in $(V(Q_\infty), \op{d}_{\rm gr})$
 and in this approximating space are within Gromov-Hausdorff distance $1$ (and a similar result
 holds for the balls in $(V(Q_n), \op{d}_{\rm gr})$). The convergences 
(\ref{D2}) and (\ref{D1}) make sense, and will indeed be proved for these approximating length spaces. 

\begin{figure}[!h]
 \begin{center}
 \includegraphics[width=10.5cm]{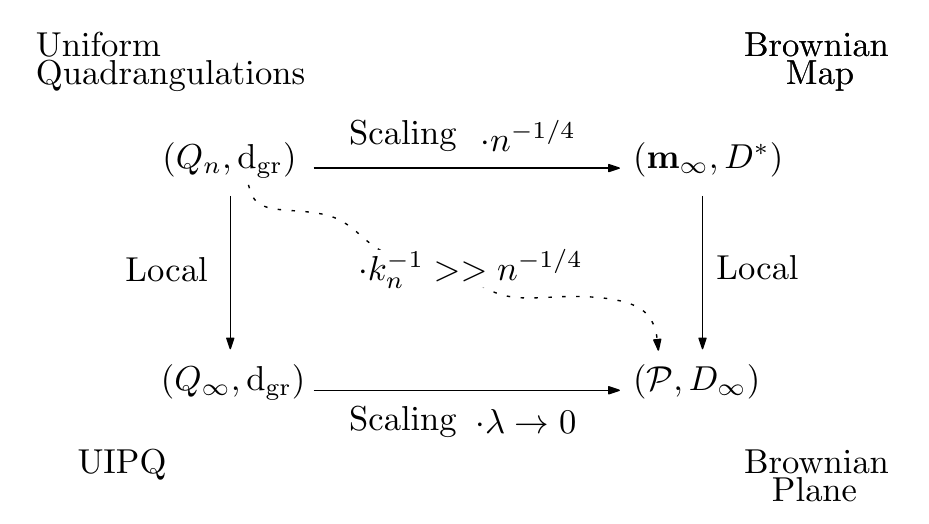}
 \caption{Illustration of Theorems  \ref{tangentcone} and \ref{scaling-limit}.}
 \end{center}
 \vspace{-4mm}
 \end{figure}

The convergences in Theorems \ref{tangentcone} and \ref{scaling-limit}, as well as the 
convergence of rescaled finite quadrangulations to the Brownian map and the local convergence to the UIPQ, 
 are summarized in the diagram of Figure 1.
 
The proofs of both Theorems \ref{tangentcone} and \ref{scaling-limit} are based on coupling methods. The coupling
argument already appears in the first statement of Theorem \ref{tangentcone}, which is
in fact much stronger than the local Gromov-Hausdorff convergence
(\ref{D0}). In the discrete setting, we prove in Proposition \ref{compmap} that we can couple $Q_{\infty}$ with $Q_{n}$ so that the balls of radius $o(n^{1/4})$ are the same with high probability. This allows us to partially answer a question of Krikun on separating cycles in the UIPQ (see Corollary \ref{Krikunconjec}).

The coupling result given in Theorem \ref{tangentcone} makes it possible to derive several important properties of the 
Brownian plane from known results about the Brownian map. In Section 5, we prove that
the Brownian plane has dimension $4$ and is homeomorphic to the plane. We also
study geodesic rays in the Brownian plane (a geodesic ray is a semi-infinite path converging to infinity whose
restriction to any finite interval is a geodesic). We prove in particular that any two geodesic rays 
coalesce in finite time, a property that is reminiscent of the results derived in \cite{AM} for the
Brownian map and in \cite{CMM} for the UIPQ. As a corollary of our study of geodesics, we also 
show that the ``labels'' $Z_a$ can be interpreted as measuring relative distances in $\pp$ from the 
point at infinity, in the sense that, a.s. for every $a,b\in \t_\infty$,
$$Z_a - Z_b = \lim_{x\to\infty} ( D_\infty(a,x)- D_\infty(b,x))$$
where the limit holds when $x$ tends to the point at infinity in the Alexandroff
compactification of $\pp$. This should be compared to Theorem 2 in \cite{CMM}. At this point it is worth mentioning that,
even if the Brownian map and the Brownian plane share many important properties, the Brownian plane enjoys the
additional scale invariance property, which may play a significant role in future investigations of these
(still mysterious) random objects. 

Finally, we conjecture that the convergence (\ref{D2}) in Theorem \ref{scaling-limit} still holds if the UIPQ $Q_\infty$
is replaced by the UIPT, and more generally that this convergence can be extended to many infinite random lattices
that are obtained as local limits of random planar maps. 

The paper is organized as follows. Section 2 gathers some preliminaries about the 
local Gromov-Hausdorff topology and the convergence of rescaled 
finite quadrangulations towards the Brownian map. Section 3 is devoted to the proof of
Theorem \ref{tangentcone}. Theorem \ref{scaling-limit} is proved in Section 4. Finally 
Section 5 studies properties of the Brownian plane. 

\medskip
\noindent{\bf Acknowledgment.} We are indebted to Gr\'egory Miermont for a number of very useful discussions.

\section{Preliminaries}

\subsection{Gromov-Hausdorff convergence}
\label{GHconv}

In this subsection, we recall some basic notions from metric geometry. For more details we refer to \cite{BBI}.

 A pointed metric space $ \mathbf{E}=(E,d,\alpha)$ is a metric space given with a distinguished point $\alpha \in E$. We will generally use bold letters $ \mathbf{E}, \mathbf{Q}, ...$  for pointed spaces. For every $r \geq 0$, we denote  the closed ball of radius $r$ centered at $\alpha$ by $B_{r}( \mathbf{E})$.  We view $B_{r}( \mathbf{E})$ as a pointed metric space (where obviously $\alpha$ is the distinguished point).
 
 Recall that if $X$ and $Y$ are two compact subsets of a metric space $(E,d)$, the Hausdorff distance between $X$ and $Y$ is 
 \begin{eqnarray*}  \mathrm{d}_{\op{H}}^E(X,Y) &:=& \inf\{ \varepsilon>0 : X \subset Y^\varepsilon \mbox{ and }Y \subset X^\varepsilon\}, \end{eqnarray*}
where $X^\varepsilon := \{ x \in E : d(x,X) \leq \varepsilon\}$ denotes the $\varepsilon$-neighborhood of $X$.
If $\mathbf{E}=(E,d,\alpha)$ and $\mathbf{E'}=(E',d',\alpha')$ are two pointed compact  metric spaces, the Gromov-Hausdorff distance between
$\mathbf{E}$ and $\mathbf{E'}$  is 
 \begin{eqnarray*} \op{d_{GH}}(\mathbf{E},\mathbf{E'}) &:=& \inf\big\{\op{d}_{\op{H}}^F(\phi(E),\phi'(E')) \vee \delta(\phi(\alpha),\phi'(\alpha')) \big\}, \end{eqnarray*} where the infimum is taken over all choices of the metric space $(F,\delta)$ and of the  isometric embeddings $\phi : E \to F$ and $\phi'  : E' \to F$ of $E$ and $E'$ into $F$. 
 The  Gromov-Hausdorff  distance is indeed a metric on the space $ \mathbb{K}$ of all isometry classes of pointed compact metric spaces (see  \cite[Section 7.4]{BBI} for more details), and the space $(\mathbb{K}, \op{d_{GH}})$ is a Polish space, that is, a separable and complete 
 metric space. 
 
 In order to extend the Gromov-Hausdorff convergence to the non-compact case, we will restrict our attention to boundedly compact length spaces.
 This restriction is not really necessary (see \cite[Section 8.1]{BBI}) but it will avoid certain technicalities, which are not relevant to
 our purposes. 
 
 If $(E,d)$ is a metric space and 
$\gamma:[0,T] \longrightarrow E$ is a continuous path in $E$, the length of $\gamma$
is defined by
  \begin{eqnarray*}L(\gamma)&:=&\sup_{0=t_{0}< \cdots<t_{k}=T}\sum_{i=0}^{k-1}d\big(\gamma(t_{i}),\gamma(t_{i+1})\big), \end{eqnarray*}
where the supremum is over all choices of the subdivision $0=t_0<t_1<\cdots<t_k=T$ of $[0,T]$. 
The space $(E,d)$ is called a length space if, for every $x,y \in E$, the distance $d(x,y)$ coincides with the infimum of the lengths of continuous paths connecting $x$ to $y$. We say that $(E,d)$ is boundedly compact if the closed balls in $(E,d)$ are compact. For a length space, this is
equivalent to saying that $(E,d)$ is complete and locally compact  \cite[Proposition 2.5.22]{BBI}.
In a boundedly compact length space, a path with minimal length (or geodesic) exists between any pair of
points.

Let $( \mathbf{E}_{n})_{n\geq0}$ and $\mathbf{E}$ be pointed boundedly compact length spaces. We say that $ \mathbf{E}_{n}$ converges
towards $ \mathbf{E}$ in the local Gromov-Hausdorff sense  if, for every $r \geq 0$, we have
  \begin{eqnarray*} \op{d_{GH}}(B_{r}( \mathbf{E}_{n}),B_{r}(\mathbf{E}))\build{\la}_{n\to\infty}^{}   0.\end{eqnarray*} 
This notion of convergence is compatible with the  distance   
\begin{eqnarray*} \mathrm{d_{LGH}}( \mathbf{E}_{1}, \mathbf{E}_{2}) &=& \sum_{k = 1}^\infty 2^{-k}\left( \mathrm{d_{GH}}(B_{k}(\mathbf{E}_{1}), B_{k}( \mathbf{E}_{2}))\wedge 1\right). \end{eqnarray*}  
As observed in \cite[Proposition 3.3]{DW}, the set $\mathbb{K}_{bcl}$ of all isometry classes of pointed boundedly compact length spaces endowed with $ \mathrm{d_{LGH}}$ is a Polish space.   

We will use the following easy consequence of the preceding considerations. If $( \mathbf{E}_{n})_{n\geq0}$ is now a sequence 
of random variables with values in $\mathbb{K}_{bcl}$, this sequence converges in distribution to $\mathbf{E}\in\mathbb{K}_{bcl}$ if and only if,
for every $r\geq 0$, $B_r(\mathbf{E}_{n})$ converges in distribution to $B_r(\mathbf{E})$ in $\mathbb{K}$. Furthermore, we may 
restrict our attention to integer values of $r$. 
\medskip

\noindent{\bf Notation.} If $ \mathbf{E}= (E,d,\alpha)$ is a pointed metric space and $\lambda$ is a positive real number, the notation $ \lambda \cdot  \mathbf{E}$
stands for the rescaled metric space $(E, \lambda\cdot d, \alpha)$.  In particular for any $\lambda, \eta >0$ we have $ \lambda \cdot B_{\eta}( \mathbf{E}) = B_{\lambda \eta}( \lambda \cdot \mathbf{E})$.

\subsection{Convergence to the Brownian map}
\label{convBromap}

In this subsection, we recall the definition of the Brownian map,
and briefly discuss the convergence of rescaled finite quadrangulations to this random metric space.
The construction of the Brownian map is very similar to that of the Brownian plane
given in the introduction, but the role of the process $X$ is now played by
a Brownian excursion.

We fix $T>0$. For reasons that will become clear later, it is convenient to index the
processes that we will define by the
parameter $\lambda=T^{1/4}$. 
Let $(\be^\lambda_t)_{0\leq t\leq T}$ be a Brownian excursion
with duration $T=\lambda^4$. For the purposes of this subsection, it would be sufficient to
take $T=\lambda=1$, but later we will deal with scaled versions of the Brownian map for 
which arbitrary values of $T$ will be useful. 

 For every $s,t\in[0,T]$, we set
$$d_{\be^\lambda}(s,t)= \be^\lambda_s + \be^\lambda_t -2 \min_{s\wedge t\leq r\leq s\vee t} \be^\lambda_r$$
and we set $s\sim_{\be^\lambda} t$ if and only if $d_{\be^\lambda}(s,t)=0$. The tree coded by
$\be^\lambda$ is the quotient space $\t_{\be^\lambda} := [0,T]/\sim_{\be^\lambda}$, which is
equipped with the distance induced by $d_{\be^\lambda}$. Note that $\t_{\be^\lambda}$
 is a scaled version of the CRT: We refer to Section 3 of \cite{LGM} for basic facts about the
 CRT and the coding of real trees by functions.
We write $p_{\be^\lambda}:[0,T]\la \t_{\be^\lambda}$ for the canonical projection, and we set 
$\rho_\lambda=p_{\be^{\lambda}}(0)=p_{\be^{\lambda}}(T)$. To simplify notation, we write $\rho=\rho_1$. 

Conditionally
given $\be^\lambda$, let $Z^\lambda=(Z^\lambda_s)_{0\leq s\leq T}$ be the centered Gaussian
process with covariance
$$E[Z^\lambda_sZ^\lambda_t\,|\, \be^\lambda]= \min_{r\in[s\wedge t,s\vee t]} \be^\lambda_r.$$
The process $(Z^\lambda_s)_{0\leq s\leq T}$ has a continuous modification, which we consider
from now on. Then, a.s. for every $s\in[0,T]$, $Z^\lambda_s$ only depends on $p_{\be^\lambda}(s)$ and so,
for every $a\in \t_{\be^\lambda}$, we may set
 $Z^\lambda_a=Z^\lambda_s$, where $s$ is an arbitrary element of $[0,T]$
such that $p_{\be^\lambda}(s)=a$.

For every $s,t\in[0,T]$ such that $s\leq t$, we set 
$$D^\circ_\lambda(s,t)= D^\circ_\lambda(t,s)=
Z^\lambda_s + Z^\lambda_t -2 \max\Big(\min_{r\in[s,t]} Z^\lambda_r , \min_{r\in [t,T]\cup[0,s]} Z^\lambda_r\Big).$$
Furthermore, we set for $a,b\in \t_{\be^\lambda}$,
$$D^\circ_\lambda(a,b) = \min\{D^\circ_\lambda (s,t) : s,t\in [0,T], \,p_{\be^\lambda}(s)=a,p_{\be^\lambda}(t)=b\},$$
and 
$$D^*_\lambda(a,b) = \inf_{a_0=a,a_1,\ldots,a_p=b} \sum_{i=1}^p D^\circ_\lambda(a_{i-1},a_i),$$
where the infimum is over all choices of the integer $p\geq 1$ and of the
finite sequence $a_0,a_1,\ldots,a_p$ in $\t_{\be ^\lambda}$ such that $a_0=a$ and
$a_p=b$.

Since clearly $D^\circ_\lambda(a,b)\geq |Z^\lambda_a-Z^\lambda_b|$, it is immediate that we have also
$$D^*_\lambda(a,b)\geq |Z^\lambda_a-Z^\lambda_b|$$
for every $a,b\in \t_{\be^\lambda}$.
More precisely, the pseudo-distance $D^*_\lambda$ satisfies the so-called ``cactus bound''
\begin{equation}
\label{cabo}
D^*_\lambda(a,b) \geq Z^\lambda_a + Z^\lambda_b -2\, \min_{c\in \llbracket a,b \rrbracket} Z^\lambda_c
\end{equation}
where $\llbracket a,b \rrbracket$ denotes the geodesic segment between
$a$ and $b$ in $\t_{\be^\lambda}$. To see this, first note that the cactus bound holds for
$D^\circ_\lambda(a,b)$ instead of $D^*_\lambda(a,b)$: This is an immediate
consequence of the definition of $D^\circ_\lambda$ and the fact that, if $p_{\be^\lambda}(s)=a$
and $p_{\be^\lambda}(t)=b$, the set $p_{\be^\lambda}([s\wedge t,s\vee t])$ must contain $\llbracket a,b \rrbracket$.
Then, given a finite sequence $a_0=a,a_1,\ldots,a_p=b$ we observe that
$\llbracket a,b \rrbracket$ is contained in the union of the segments 
$\llbracket a_{i-1},a_i \rrbracket$ for $1\leq i\leq p$. Hence, there must exist
an index $j\in\{1,\ldots,p\}$ such that
$$\min_{c\in \llbracket a,b \rrbracket} Z^\lambda_c \geq \min_{c\in \llbracket a_{j-1},a_j \rrbracket} Z^\lambda_c$$
and it follows that
\begin{align*}
\sum_{i=1}^p D^\circ_\lambda(a_{i-1},a_i) &\geq \sum_{i=1}^{j-1} |Z^\lambda_{a_i} - Z^\lambda_{a_{i-1}}| +
( Z^\lambda_{a_{j-1}} + Z^\lambda_{a_j} - 2\min_{c\in \llbracket a_{j-1},a_j \rrbracket} Z^\lambda_c )
+  \sum_{i=j+1}^{p} |Z^\lambda_{a_i} - Z^\lambda_{a_{i-1}}|\\
&\geq Z^\lambda_a + Z^\lambda_b -2\, \min_{c\in \llbracket a,b \rrbracket} Z^\lambda_c.
\end{align*}

Finally, we put $a\approx_\lambda b$ if and only if $D^*_\lambda(a,b)=0$. Specializing now to the case $T=1$,  the Brownian map $\mathrm{m}_\infty$
is defined as the quotient space $\t_{\be^1}/\approx_1$, which is equipped with the distance induced by
$D^*_1$. We view $\mathrm{m}_\infty$ as a pointed metric space with distinguished 
point $\rho$ (here and later we abuse notation and identify $\rho$ with its equivalence class in $\t_{\be^1}/\approx_1$), and we write
$\mathbf{m}_\infty=(\mathrm{m}_\infty, D^*_1,\rho)$ for this pointed metric space. With the notation 
introduced at the end of the preceding subsection, we have, by a simple scaling argument,
\begin{equation}
\label{identdist}
\lambda\cdot \mathbf{m}_\infty \build{=}_{}^{\rm(d)} \left( \t_{\be^\lambda}/\approx_\lambda, D^*_\lambda, \rho_\lambda\right),
\end{equation}
with the same abuse of notation for $\rho_\lambda$. This identity in distribution explains why
we considered an arbitrary value of $T>0$ (and not only the case $T=1$) in the preceding discussion.

 Let us recall the convergence of rescaled random quadrangulations towards the Brownian map \cite{LGU,Mi-quad}. For every integer 
 $n\geq 1$, let $Q_n$ and  $V(Q_n)$ be as in the introduction. We view $V(Q_n)$ as a metric space for the graph distance
 $\op{d}_{\rm gr}$, which is pointed at the root vertex $\rho_{(n)}$ of $Q_n$. 
Then, we have
 \begin{equation}
 \label{cvbm}
(V(Q_n),n^{-1/4}\op{d}_{\rm gr}, \rho_{(n)}) \xrightarrow[n\to\infty]{(d)} \Big(\mathrm{m}_\infty,(\frac{8}{9})^{1/4}D^*_1, \rho\Big)
 \end{equation}
  in distribution in $\mathbb{K}$.  

\section{The Brownian plane as a tangent cone of the Brownian map}
\label{secasympcone}

In this section, we prove Theorem \ref{tangentcone}. To do so, we establish a coupling result (Proposition \ref{couplingBMBP}) from which Theorem  \ref{tangentcone} easily follows. This lemma shows that we can couple the realizations of the Brownian plane  and of the Brownian map in such a way that 
the balls of small radius centered at the distinguished point coincide with high probability. 

\subsection{Absolute continuity properties of excursion measures}
\label{absoco}

The branching structures of the (scaled) CRT $\t_{\be^\lambda}$ and of the infinite 
Brownian tree $\t_\infty$ are encoded respectively in a Brownian excursion $\be^\lambda$
of duration $T=\lambda^4$ and in a pair $(R,R')$ of independent three-dimensional Bessel processes. In the proof 
of the forthcoming coupling result, it will be important to have information about the
absolute continuity properties of the law of a pair of processes corresponding respectively to the initial and
the final part of $\be^\lambda$, with respect to the law
of $(R,R')$.

We fix $T>0$ and set $\lambda=T^{1/4}$ as previously. We also consider two
reals $\alpha,\beta>0$ such that $\alpha + \beta <T$. We will then consider the
pair $((\be^\lambda_t)_{0\leq t\leq \alpha}, (\be^\lambda_{T-t})_{0\leq t\leq \beta})$, which is viewed
as a random element of the space $C([0,\alpha],\R_+)\times C([0,\beta],\R_+)$. The generic
element of the latter space will be denoted by $(\omega,\omega')$.

For $t>0$ and $y\in \R$, we let 
$$p_t(y)= \frac{1}{\sqrt{2\pi t}} \,\exp(-\frac{y^2}{2t})$$
be the usual Brownian transition density. The transition density 
of Brownian motion killed upon hitting $0$ is 
$$p^*_t(x,y)= p_t(y-x) - p_t(y+x)$$
for $x,y>0$. We also set
$$q_t(x) = \frac{x}{t} p_t(x)$$
for $t>0$ and $x>0$. We recall that the transition density of
the three-dimensional Bessel process is given by
$$\left\{
\begin{array}{ll}
p^{(3)}_t(x,y) = {\displaystyle\frac{y}{x}}\,p^*_t(x,y)\qquad &\hbox{if }t>0, x>0, y\geq 0,\\
\noalign{\smallskip}
p^{(3)}_t(0,y) = 2y\,q_t(y)&\hbox{if }t>0, y\geq 0.
\end{array}
\right.
$$

\begin{proposition}
\label{absocont}
The law of the triplet
$$\Big( (\be^\lambda_{t})_{0\leq t\leq \alpha}, (\be^\lambda_{T-t})_{0\leq t\leq \beta},
\min_{\alpha\leq t\leq T-\beta} \be^\lambda_t\Big)$$
is absolutely continuous with respect to the law of
$$\Big( (R_t)_{0\leq t\leq \alpha}, (R'_t)_{0\leq t\leq \beta} ,
\inf_{t\geq \alpha} R_t \wedge \inf_{t\geq \beta} R'_t\Big),$$
with density given by the function $\Delta_{T,\alpha,\beta}(\omega(\alpha),\omega'(\beta),z)$,
where, for every $x,y,z>0$,
$$\Delta_{T,\alpha,\beta}(x,y,z)= \varphi_{T,\alpha,\beta}(x,y)\,
\psi_{T-(\alpha +\beta)} (x,y,z)$$
with
$$\varphi_{T,\alpha,\beta}(x,y) = \frac{\sqrt{2\pi T^3}}{2xy} \,p^*_{T-(\alpha+\beta)}(x,y)$$
and, for every $s>0$, 
$$\psi_s(x,y,z) = \frac{
(2(x+y-2z)/s)\,\exp(-((x+y-2z)^2 -(y-x)^2)/2s)}
{(1-\exp(-2xy/s))\,(\frac{1}{x} + \frac{1}{y} -\frac{2z}{xy})}
\;{\bf 1}_{\{z<x\wedge y\}}.$$
Moreover, for every $x,y,z>0$ such that $z<x\wedge y$, we have
$$\lim_{T\to\infty} \Delta_{T,\alpha,\beta}(x,y,z) = 1.$$
\end{proposition}

\proof We first recall a well-known fact about the Brownian excursion
with fixed duration. 
If $0<t_1<t_2<\cdots< t_p<T$, the density of the law of the $p$-tuple 
$(\be^\lambda_{t_1},\ldots,\be^\lambda_{t_p})$ is
$$2\sqrt{2\pi T^3}\, q_{t_1}(x_1)p^*_{t_2-t_1}(x_1,x_2)\ldots p^*_{t_p-t_{p-1}}(x_{p-1},x_p)
q_{T-t_p}(x_p).$$
See e.g. Chapter XII of \cite{RY} for a proof. To simplify notation, we set,
for $0<t_1<\cdots <t_p$ and $x_1,\ldots,x_p>0$,
$$F_{t_1,\ldots,t_p}(x_1,\ldots,x_p)= q_{t_1}(x_1)p^*_{t_2-t_1}(x_1,x_2)\ldots p^*_{t_p-t_{p-1}}(x_{p-1},x_p).$$
Then fix $0<t_1<\cdots <t_p=\alpha$ and $0<t'_1<\cdots <t'_q=\beta$. We see that the density of
$(\be^\lambda_{t_1},\ldots,\be^\lambda_{t_p},\be^\lambda_{T-t'_q},\ldots,\be^\lambda_{T-t'_1})$ is the function
\begin{align*}
&g^T_{t_1,\ldots,t_p,t'_1,\ldots,t'_q}(x_1,\ldots,x_{p+q})\\
&\qquad= 2\sqrt{2\pi T^3}\,  F_{t_1,\ldots,t_p}(x_1,\ldots,x_p)
\,F_{t'_1,\ldots,t'_q}(x_{p+q},\ldots,x_{p+1})\,p^*_{T-(\alpha+\beta)}(x_p,x_{p+1}).
\end{align*}
On the other hand, from the explicit formulas for the transition density of the
three-dimensional Bessel process, we have 
\begin{align*}
F_{t_1,\ldots,t_p}(x_1,\ldots,x_p)&= \frac{1}{2x_p} \,
p^{(3)}_{t_1}(0,x_1)p^{(3)}_{t_2-t_1}(x_1,x_2)\cdots p^{(3)}_{t_p-t_{p-1}}(x_{p-1},x_p)
= \frac{1}{2x_p} \, G^{(3)}_{t_1,\ldots,t_p}(x_1,\ldots,x_p),
\end{align*}
where $G^{(3)}_{t_1,\ldots,t_p}$ is the density of $(R_{t_1},\ldots,R_{t_p})$. Consequently, we have also
\begin{align*}
&g^T_{t_1,\ldots,t_p,t'_1,\ldots,t'_q}(x_1,\ldots,x_{p+q})\\
&\qquad= 2\sqrt{2\pi T^3}\, G^{(3)}_{t_1,\ldots,t_p}(x_1,\ldots,x_p)\,G^{(3)}_{t'_1,\ldots,t'_q}(x_{p+q},\ldots,x_{p+1})\times 
\frac{p^*_{T-(\alpha+\beta)}(x_p,x_{p+1})}{4x_px_{p+1}}.
\end{align*}
It follows from this calculation that the density of the law of the pair
$$\Big( (\be^\lambda_{t})_{0\leq t\leq \alpha}, (\be^\lambda_{T-t})_{0\leq t\leq \beta}\Big)$$
with respect to the law of
$$\Big( (R_t)_{0\leq t\leq \alpha}, (R'_t)_{0\leq t\leq \beta}\Big),$$
is given by the function $\varphi_{T,\alpha,\beta}(\omega(\alpha),\omega'(\beta))$.

We next consider the conditional density of the third component 
of each triplet in the proposition, given its first two components. We note that,
for every $z>0$, 
\begin{align*}
P\Big[ \inf_{t\geq \alpha} R_t \wedge \inf_{t\geq \beta} R'_t > z\;\Big|\; (R_t)_{0\leq t\leq \alpha}, (R'_t)_{0\leq t\leq \beta}\Big]
&= P\Big[ \inf_{t\geq \alpha} R_t > z \;\Big|\; R_\alpha\Big] \; P\Big[ \inf_{t\geq \beta} R'_t > z \;\Big|\; R'_\beta\Big] \\
&= \Big( 1 - \frac{z\wedge R_\alpha}{R_\alpha}\Big) \Big( 1 - \frac{z\wedge R'_\beta}{R'_\beta}\Big) .
\end{align*}
Hence the conditional density of 
$$ \inf_{t\geq \alpha} R_t \wedge \inf_{t\geq \beta} R'_t $$
given $(R_t)_{0\leq t\leq \alpha}$ and $ (R'_t)_{0\leq t\leq \beta}$ is
$$h_{R_\alpha,R'_\beta}(z)= {\bf 1}_{\{0<z<R_\alpha\wedge R'_\beta\}} \, 
\Big(\frac{1}{R_\alpha} + \frac{1}{R'_{\beta}} - \frac{2z}{R_\alpha R'_\beta}\Big).
$$

Similarly, to get the conditional distribution of 
$$\min_{\alpha\leq t\leq T-\beta} \be^\lambda_t,$$
we observe that, conditionally on $(\be^\lambda_t)_{0\leq t\leq \alpha}$
and $(\be^\lambda_{T-t})_{0\leq t\leq \beta}$, the process
$(\be^\lambda_{\alpha+t})_{0\leq t\leq T-(\alpha+\beta)}$ is a 
Brownian bridge of duration $T-(\alpha+\beta)$, starting from $\be^\lambda_\alpha$
and ending at $\be^\lambda_{T-\beta}$, and conditioned not to hit $0$. 

Fix $s>0$ and $x>0$, and recall that for a standard Brownian motion $B$
starting from $x$ the density of the pair $(B_s,\min_{0\leq r\leq s} B_r)$
is the function $g_s(y,z)=2q_s(x+y-2z)\,{\bf 1}_{\{z<x\wedge y\}}$. Hence,
if we also fix $y>0$, the conditional 
density of $\min_{0\leq r\leq s} B_r$ knowing that $B_s=y$ is 
$$f_{s,x,y}(z) = \frac{ 2q_s(x+y-2z)}{ p_s(y-x)}\,{\bf 1}_{\{z<x\wedge y\}}.$$
In particular,
$$P\Big(\min_{0\leq r\leq s} B_r >0\;\Big|\; B_s=y\Big)
= \int_0^{x\wedge y} dz\,f_{s,x,y}(z)  = 1-\exp(-\frac{2xy}{s}).$$
We conclude from these calculations that the density of $\min_{\alpha\leq t\leq T-\beta} \be^\lambda_t$
knowing $(\be^\lambda_t)_{0\leq t\leq \alpha}$
and $(\be^\lambda_{T-t})_{0\leq t\leq \beta}$ is the function $f^*_{T-(\alpha+\beta),\be^\lambda_\alpha,\be^\lambda_{T-\beta}}(z)$,
where, for $s>0$ and $x,y,z>0$,
$$f^*_{s,x,y}(z) = \frac{f_{t,x,y}(z)}{1-\exp(-\frac{2xy}{s})}\,{\bf 1}_{\{0<z<x\wedge y\}}.$$

Set 
$$m_{\alpha,\beta}=\inf_{t\geq \alpha} R_t \wedge \inf_{t\geq \beta} R'_t$$
to simplify notation. It follows that, for every nonnegative measurable 
function $\Gamma$ on $C([0,\alpha],\R_+)\times C([0,\beta],\R_+)\times\R$,
\begin{align*}
&E\Big[\Gamma\Big((\be^\lambda_{t})_{0\leq t\leq \alpha}, (\be^\lambda_{T-t})_{0\leq t\leq \beta},
\min_{\alpha\leq t\leq T-\beta} \be^\lambda_t\Big)\Big]\\
&=E\Big[ \int_0^{\infty}
dz\,f^*_{T-(\alpha+\beta),\be^\lambda_\alpha,\be^\lambda_{T-\beta}}(z)\,
\Gamma\Big((\be^\lambda_{t})_{0\leq t\leq \alpha}, (\be^\lambda_{T-t})_{0\leq t\leq \beta},z\Big)\Big]\\
&=E\Big[\varphi_{T,\alpha,\beta}(R_\alpha,R'_\beta)
\int_0^{\infty}
dz\,f^*_{T-(\alpha+\beta),R_\alpha,R'_{\beta}}(z)\,
\Gamma\Big((R_t)_{0\leq t\leq \alpha}, (R'_t)_{0\leq t\leq \beta},z\Big)\Big]\\
&= E\Big[ \int_0^\infty dz \,h_{R_\alpha,R'_\beta}(z)\,\Delta_{T,\alpha,\beta}(R_\alpha,R'_\beta,z)\,
\Gamma\Big((R_t)_{0\leq t\leq \alpha}, (R'_t)_{0\leq t\leq \beta},z\Big)\Big]\\
&= E\Big[ \Delta_{T,\alpha,\beta}(R_\alpha,R'_\beta,m_{\alpha,\beta})\,
\Gamma\Big((R_t)_{0\leq t\leq \alpha}, (R'_t)_{0\leq t\leq \beta},m_{\alpha,\beta}\Big)\Big].
\end{align*}
In the second equality, we applied the first part of the proof, and in the third one, we used the
identity
$$\varphi_{T,\alpha,\beta}(R_\alpha,R'_\beta)f^*_{T-(\alpha+\beta),R_\alpha,R'_{\beta}}(z)
= h_{R_\alpha,R'_\beta}(z)\,\Delta_{T,\alpha,\beta}(R_\alpha,R'_\beta,z),$$
which follows from our definitions. The proof of the first assertion of the proposition is now complete. 
The last assertion follows from the explicit expressions for $\varphi_{T,\alpha,\beta}(x,y)$ and 
$\psi_{T-(\alpha +\beta)} (x,y,z)$. \endproof

\subsection{Proof of Theorem \ref{tangentcone}}

In this section, we prove Theorem \ref{tangentcone}. 
Recall the notation $\boldsymbol{ \mathcal{P}}=(\pp, D_\infty,\rho_\infty)$ for the Brownian plane
viewed as a pointed metric space.

\begin{proposition}
\label{couplingBMBP}
Let $\delta>0$ and $r\geq 0$. There exists $\lambda_0>0$ such that, for every $\lambda\geq \lambda_0$, we  can construct copies of $\lambda \cdot \mathbf{m}_\infty$
and $\boldsymbol{ \mathcal{P}}$ on the same probability space, in such a way that the equality
\begin{equation}
\label{idenball}
B_r(\lambda \cdot \mathbf{m}_\infty) = B_r(\boldsymbol{ \mathcal{P}})
\end{equation}
holds with probability at least $1-\delta$. 
\end{proposition}

Proposition \ref{couplingBMBP} immediately implies that, for every $r\geq 0$,
\begin{equation}
\label{converball}
B_r(\lambda\cdot \mathbf{m}_\infty)\build{\la}_{\lambda\to\infty}^{(d)} B_r(\boldsymbol{ \mathcal{P}})
\end{equation}
in distribution in $\mathbb{K}$. From the observations of the end of subsection \ref{GHconv}, this suffices to
get the convergence (\ref{D0}) in Theorem \ref{tangentcone}. The first assertion in Theorem \ref{tangentcone} also
readily follows from Proposition \ref{couplingBMBP} and the scale invariance of the Brownian plane.

\medskip
\proof
We will
rely on the identity in distribution \eqref{identdist}.
To simplify 
notation, we write $Y^\lambda=\t_{\be^{\lambda}}/\approx_{\lambda}$. 
We also write $\mathbf{Y}^\lambda$ for the pointed space
$(Y^\lambda, D^*_{\lambda}, \rho_\lambda)$, so that $\lambda\cdot \mathbf{m}_\infty$ has the same distribution
as $\mathbf{Y}^\lambda$. 

 We first choose $A>1$ sufficiently large so that, if
$\beta=(\beta_t)_{t\geq 0}$ is a standard linear Brownian motion, we have
$$P\Big[\Big\{\min_{0\leq t\leq A} \beta_t < -10r\Big\}\cap \Big\{\min_{A\leq t\leq A^2} \beta_t < -10r\Big\}
\cap \Big\{\min_{A^2\leq t\leq A^4} \beta_t < -10r\Big\}\Big] \geq 1-\delta/8.$$
Next we choose $\alpha >0$ large enough so that 
the property
$$\inf_{t\geq \alpha} R_t \wedge \inf_{t\geq \alpha} R'_t > A^4$$
holds with probability at least $1-\delta/8$.

From  Proposition \ref{absocont} (especially the last assertion of this proposition) and  standard coupling arguments, we can 
find $\lambda_0>0$, such that $\lambda_0^4 > 2\alpha$, and such that, for every
$\lambda\geq \lambda_0$, we can construct both processes $\be^{\lambda}$ and $(R,R')$ 
on the same probability space, in such a way that the probability of the
event
$$\mathcal{E}_\lambda:=\{\be^{\lambda}_t=R_t\hbox{ and } \be^{\lambda}_{\lambda^4-t}=R'_t, \forall t\leq \alpha\}
\cap\Big\{\min_{\alpha\leq t\leq \lambda^4-\alpha} \be^{\lambda}_t = \inf_{t\geq \alpha} R_t \wedge \inf_{t\geq \alpha} R'_t\Big\}$$
is bounded below by $1-\delta/2$. 

We fix $\lambda\geq \lambda_0$
and we write $T=\lambda^4$ as in subsection \ref{absoco}. We suppose that $\be^\lambda$ and 
the pair $(R,R')$ have been constructed so that $P(\mathcal{E}_\lambda)>1-\delta/2$. 
We then observe that, conditionally given $\be^{\lambda}$, the process $(W^{\lambda}_t)_{t\in[-\alpha,\alpha]}$
given by
$$ W^{\lambda}_t = \left\{ \begin{array}{ll} Z^{\lambda}_t\qquad&\hbox{if }t\in[0,\alpha]\\
Z^{\lambda}_{T+t}\qquad&\hbox{if }t\in[-\alpha,0]
\end{array}
\right.$$
is a Gaussian process whose covariance is a function of the triplet
$$\Big( (\be^{\lambda}_{t})_{0\leq t\leq \alpha}, (\be^{\lambda}_{T-t})_{0\leq t\leq \alpha},
\min_{\alpha\leq t\leq T-\alpha} \be^{\lambda}_t\Big).$$
Similarly, conditionally given $(R,R')$, the process $(Z_t)_{t\in[-\alpha,\alpha]}$
is also Gaussian with covariance given by the {\em same} function of the triplet
$$\Big( (R_t)_{0\leq t\leq \alpha}, (R'_t)_{0\leq t\leq \alpha} ,
\inf_{t\geq \alpha} R_t \wedge \inf_{t\geq \alpha} R'_t\Big).$$
Note that $\mathcal{E}_\lambda$ is precisely the event where the latter two triplets of processes
coincide. It follows that we can construct the processes $Z^{\lambda}$ and $Z$ 
in such a way that, on the event $\mathcal{E}_\lambda$, we have also the identities
$$Z^{\lambda}_t = Z_t\quad,\quad Z^{\lambda}_{T-t}= Z_{-t}\quad,\ \forall t\in[0,\alpha].$$
From now on, we assume that these identities hold on $\mathcal{E}_\lambda$.

To simplify notation, we write $\t_{(\lambda)}= \t_{\be^{\lambda}}$,  and $p_{(\lambda)}$
for the canonical projection from $[0,T]$ onto $\t_{(\lambda)}$. For every $x\in[0,\be^{\lambda}_{T/2}]$, we set
$$\gamma_\lambda(x)= \sup\{t\leq T/2 : \be^{\lambda}_t=x\}\quad,\quad \eta_\lambda(x)= \inf\{t\geq T/2:\be^{\lambda}_t=x\}.$$
Then $p_{(\lambda)}(\gamma_\lambda(x))=p_{(\lambda)}(\eta_\lambda(x))$ is the vertex at distance $x$ from the root in the ancestral
line of the vertex $p_{(\lambda)}(T/2)$ in the tree $\t_{(\lambda)}$. Similarly, we define, for every $x\geq 0$,
$$\gamma_\infty(x)=\sup\{t\geq 0: R_t=x\}\quad,\quad \eta_\infty(x)=\sup\{t\geq 0: R'_t=x\}.$$
With the notation of the introduction, we have $p_\infty(\gamma_\infty(x))=p_\infty(-\eta_\infty(x))$
for every $x\geq 0$. Furthermore, the process $(Z_{\gamma_\infty(x)})_{x\geq 0}$ is distributed as 
a standard linear Brownian motion. 

Write $\mathcal{F}_\lambda$ for the intersection of $\mathcal{E}_\lambda$ with the event where we have both
$$\inf_{t\geq \alpha} R_t \wedge \inf_{t\geq \alpha} R'_t > A^4$$
and 
\begin{equation}
\label{tech1}
\min_{0\leq x\leq A} Z_{\gamma_\infty(x)} < -10r\quad,\quad\min_{A\leq x\leq A^2} Z_{\gamma_\infty(x)}< -10r\quad,\quad
\min_{A^2\leq x\leq A^4} Z_{\gamma_\infty(x)}< -10r.
\end{equation}
By our choice of the constants $A$ and $\alpha$, we have $P(\mathcal{F}_\lambda)\geq 1-\delta$. 

On the event $\mathcal{F}_\lambda$, we have
$$\min_{\alpha\leq t\leq T-\alpha} \be^\lambda_t = \inf_{t\geq \alpha} R_t \wedge \inf_{t\geq \alpha} R'_t>A^4$$
and therefore, for every $x\in[0,A^4]$,
$\gamma_\lambda(x)=\gamma_\infty(x)<\alpha$
and 
$T-\eta_\lambda(x)= \eta_\infty(x) <\alpha$.
It follows that, for every $x\in [0,A^4]$,
$$Z^\lambda_{\gamma_\lambda(x)}= Z_{\gamma_\infty(x)}=Z_{-\eta_\infty(x)}= Z^\lambda_{\eta_\lambda(x)}.$$

Before stating another lemma, we introduce one more piece of notation. Let
$s,t\in[0,T]$.
If $s$ and $t$ both belong to $[0,T/2]$ or if they both belong to $[T/2,T]$, we set
$$\wt D^\circ_{\lambda}(s,t) = Z^\lambda_s + Z^\lambda_t - 2 \min_{u\in [s\wedge t,s\vee t]} Z^\lambda_u.$$
Otherwise we set
$$\wt D^\circ_{\lambda}(s,t) = Z^\lambda_s + Z^\lambda_t - 2 \min_{u\in [0,s\wedge t]\cup [s\vee t,T]} Z^\lambda_u.$$

\begin{lemma}
\label{keyasymp}
Suppose that $\mathcal{F}_\lambda$ holds. 

\noindent {\rm (i)} For every $t\in [\gamma_\lambda(A),\eta_\lambda(A)]$,
$$D^*_{\lambda}(\rho_{\lambda},p_{(\lambda)}(t)) > r.$$
For every $s,t\in [0,\gamma_\lambda(A)]\cup [\eta_\lambda(A),T]$ such that
$D^*_{\lambda}(\rho_{\lambda},p_{(\lambda)}(s))\leq r$ and $D^*_{\lambda}(\rho_{\lambda},p_{(\lambda)}(t))\leq r$,
we have
\begin{equation}
\label{tech10}
D^*_{\lambda}(p_{(\lambda)}(s),p_{(\lambda)}(t)) = \inf_{s_0,t_0,s_1,t_1,\ldots,s_p,t_p} 
\sum_{i=1}^p \wt D^\circ_{\lambda}(t_{i-1},s_i)
\end{equation}
where the infimum is over all choices of the integer $p\geq 0$ and of the
reals $s_0,s_1,\ldots,s_p,t_0,t_1,\ldots,t_{p}$ belonging to $[0,\gamma_\lambda(A^2)]\cup [\eta_\lambda(A^2),T]$ such that 
$s_0=s$, $t_p=t$, and $p_{(\lambda)}(s_i)=p_{(\lambda)}(t_i)$ for every 
$i\in\{0,1,\ldots,p\}$.

\noindent{\rm (ii)} For every $t'\in (-\infty, -\eta_\infty(A)] \cup [\gamma_\infty(A),\infty)$,
$$D_\infty(\rho_\infty,p_\infty(t')) > r.$$
For every $s',t'\in [-\eta_\infty(A),\gamma_\infty(A)]$ such that
$D_\infty(\rho_\infty,p_\infty(s')) \leq r$ and $D_\infty(\rho_\infty,p_\infty(t')) \leq r$,
we have
\begin{equation}
\label{tech11}
D_\infty(p_\infty(s'),p_\infty(t'))= \inf_{s'_0,t'_0,s'_1,t'_1,\ldots,s'_p,t'_p} 
\sum_{i=1}^p D^\circ_{\infty}(t'_{i-1},s'_i)
\end{equation}
where the infimum is over all choices of the integer $p\geq 0$ and of the
reals $s'_0,s'_1,\ldots,s'_p,t'_0,t'_1,\ldots,t'_{p}$ belonging to $[-\eta_\infty(A^2),\gamma_\infty(A^2)]$ such that 
$s'_0=s'$, $t'_p=t'$, and $p_{\infty}(s'_i)=p_{\infty}(t'_i)$ for every 
$i\in\{0,1,\ldots,p\}$.
\end{lemma}

\proof 
We argue on the event $\mathcal{F}_\lambda$ and prove (i). We observe that, for every
$t\in[\gamma_\lambda(A), \eta_\lambda(A)]$, we have
$$D^*_{\lambda}(\rho_\lambda,p_{(\lambda)}(t))\geq Z^\lambda_t - 2\min_{c\in \llbracket \rho_\lambda,p_{(\lambda)}(t)\rrbracket} Z^\lambda_c
\geq - \min_{c\in \llbracket \rho_\lambda,p_{(\lambda)}(t)\rrbracket} Z^\lambda_c \geq 10r.$$
The first inequality is the cactus bound \eqref{cabo}, and the last one follows from the fact
that if $t\in[\gamma_\lambda(A), \eta_\lambda(A)]$, the ancestral line of $p_{(\lambda)}(t)$ contains all vertices
of the form $p_{(\lambda)}(\gamma_\lambda(x))$ for $x\in[0,A]$, and we use the equality
$Z^\lambda_{\gamma_\lambda(x)}= Z_{\gamma_\infty(x)}$ (which holds for these values of $x$ on $\mathcal{F}_\lambda$), together with the first bound in (\ref{tech1}). 

We next turn to the proof of the second assertion in (i). Let $s,t\in [0,\gamma_\lambda(A)]\cup [\eta_\lambda(A),T]$ such that
$D^*(\rho_{\lambda},p_{(\lambda)}(s))\leq r$ and $D^*(\rho_{\lambda},p_{(\lambda)}(t))\leq r$.
Note that we have then
$$|Z^\lambda_s| \leq r\quad,\quad |Z^\lambda_t|\leq r$$
by the cactus bound. Furthermore, we have 
by definition
\begin{equation}
\label{tech12}
D^*_{\lambda}(p_{(\lambda)}(s),p_{(\lambda)}(t)) = \inf_{s=s_0,t_0,s_1,t_1,\ldots,s_p,t_p=t} 
\sum_{i=1}^p D^\circ_{\lambda}(t_{i-1},s_i)
\end{equation}
where the infimum is over all choices of the integer $p\geq 0$ and of the
reals $s_0,s_1,\ldots,s_p,t_0,t_1,\ldots,t_p$ in $[0,T]$ such that 
$s_0=s$, $t_p=t$ and $p_{(\lambda)}(s_i)=p_{(\lambda)}(t_i)$ for every 
$i\in\{0,1,\ldots,p\}$.  

Since
$D^*_{\lambda}(p_{(\lambda)}(s),p_{(\lambda)}(t))\leq 2r$, we can obviously restrict our attention
to reals $s_0,s_1,\ldots,t_p$ such that
\begin{equation}
\label{tech5}
\sum_{i=1}^p D^\circ_{\lambda}(t_{i-1},s_i) < 5r/2.
\end{equation}
We next claim that in the infimum in the 
right-hand side of (\ref{tech12}), we can furthermore limit ourselves to
choices of $s_1,\ldots,s_p,t_0,t_1,\ldots,t_{p-1}$ such that, for every
$i\in\{0,1,\ldots,p\}$, both $s_{i}$ and $t_i$ belong to
$[0,\gamma_\lambda(A^2)]\cup [\eta_\lambda(A^2),T]$.

Indeed, suppose that, for some
$i\in\{0,1,\ldots,p\}$, $t_i$ (or $s_i$) does not belong to
$[0,\gamma_\lambda(A^2)]\cup [\eta_\lambda(A^2),T]$. From (\ref{tech5}), we have $D^*_{\lambda}(p_{(\lambda)}(s), p_{(\lambda)}(t_i))< 5r/2$. On the other
hand, by the cactus bound and the property $Z^\lambda_s\geq -r$, we have
$$D^*_{\lambda}(p_{(\lambda)}(s), p_{(\lambda)}(t_i)) \geq -r - \min_{c\in \llbracket p_{(\lambda)}(s),p_{(\lambda)}(t_i)\rrbracket} Z^\lambda_c \geq 9r$$
because the fact that $s\in [0,\gamma_\lambda(A)]\cup [\eta_\lambda(A),T]$ and 
$t_i\notin [0,\gamma_\lambda(A^2)]\cup [\eta_\lambda(A^2),T]$ ensures that the geodesic
segment $\llbracket p_{(\lambda)}(s),p_{(\lambda)}(t_i)\rrbracket$ contains all vertices
of the form $p_{(\lambda)}(\gamma_\lambda(x))$ for $x\in[A,A^2]$, and we use the equality
$Z^\lambda_{\gamma_\lambda(x)}= Z_{\gamma_\infty(x)}$
(for these values of $x$), together with the second bound in (\ref{tech1}).
This contradiction proves our claim. 

In order to establish formula (\ref{tech10}), we still need to justify the fact that
we can replace $D^\circ_{\lambda}$ by $\wt D^\circ_{\lambda}$ in (\ref{tech12}).
So let $s_0,t_1,\ldots,t_p$ be reals in $[0,\gamma_\lambda(A^2)]\cup [\eta_\lambda(A^2),T]$ such 
that $s_0=s$, $t_p=t$ and (\ref{tech5}) holds. 
Consider first the case where $i\in\{1,\ldots,p\}$ is such that $t_{i-1}\in [0,\gamma_\lambda(A^2)]$
and $s_i\in [\eta_\lambda(A^2),T]$. 
Using (\ref{tech5}) we have 
$D^*_{\lambda}(p_{(\lambda)}(s),p_{(\lambda)}(t_{i-1})) < 5r/2$, hence
$Z^\lambda_{t_{i-1}} \geq Z^\lambda_s - 5r/2 \geq -7r/2$, and similarly
$Z^\lambda_{s_i} \geq -7r/2$. We have also
$$-\min_{u\in [t_{i-1},s_i]} Z^\lambda_u
\geq 10r $$
because the interval $[t_{i-1},s_i]$ must contain all $\gamma_{\infty}(x)$ for $x\in [A^2,A^4]$, and we 
use the third bound of (\ref{tech1}). Then
\begin{align*}
5r/2 > D^\circ_{\lambda}(t_{i-1},s_i)&= Z^\lambda_{t_{i-1}} + Z^\lambda_{s_i} - 2 \max \Big(\min_{u\in [t_{i-1},s_i]} Z^\lambda_u,
\min_{u\in[0,t_{i-1}]\cup [s_i,T]} Z^\lambda_u\Big)\\
&\geq -7r - 2 \max \Big(\min_{u\in [t_{i-1},s_i]} Z^\lambda_u,
\min_{u\in[0,t_{i-1}]\cup [s_i,T]} Z^\lambda_u\Big)
\end{align*}
and the previous two displays can hold only if 
$$\max \Big(\min_{u\in [t_{i-1},s_i]} Z^\lambda_u,
\min_{u\in[0,t_{i-1}]\cup [s_i,T]} Z^\lambda_u\Big)= \min_{u\in[0,t_{i-1}]\cup [s_i,T]} Z^\lambda_u,$$
which means that $\wt D^\circ_{\lambda}(t_{i-1},s_i)= D^\circ_{\lambda}(t_{i-1},s_i)$. 
Next consider the case where 
both $t_{i-1}$ and $s_i$ belong to $[0,\gamma_\lambda(A^2)]$. Then,
by the same argument as previously, we have 
$$-\min_{u\in [t_{i-1}\vee s_i, T]} Z^\lambda_u
\geq 10r $$
and it again follows that $\wt D^\circ_{\lambda}(t_{i-1},s_i)= D^\circ_{\lambda}(t_{i-1},s_i)$.
The other cases are treated in a similar way.
This completes
the proof of assertion (i).

The proof of assertion (ii) is similar. Just note that a version of the cactus bound \eqref{cabo} 
where $D^*_\lambda$ is replaced by $D_\infty$
and $Z^\lambda$ by $Z$
holds for $a,b\in\t_\infty$,  with exactly the same proof. We omit the details. \endproof

The next lemma is a simple corollary of Lemma \ref{keyasymp}.

\begin{lemma}
\label{corokey}
Assume that $\mathcal{F}_\lambda$ holds. Let $s,t\in [0,\gamma_\lambda(A)]\cup [\eta_\lambda(A),T]$. Set
$s'=s$ if $s \in [0,\gamma_\lambda(A)]$ and $s'=s-T$ if $s\in [\eta_\lambda(A),T]$,
and similarly $t'=t$  if $t \in [0,\gamma_\lambda(A)]$ and $t'=t-T$ if $t\in [\eta_\lambda(A),T]$.
Then we have $D^*_{\lambda}(\rho_{\lambda},p_{(\lambda)}(s))\leq r$ and 
$D^*_{\lambda}(\rho_{\lambda},p_{(\lambda)}(t))\leq r$ if and only if 
$D_\infty(\rho_\infty,p_\infty(s'))\leq r$ and $D_\infty(\rho_\infty,p_\infty(t'))\leq r$.
Furthermore, if these conditions hold, we have also 
$$D^*_{\lambda}(p_{(\lambda)}(s),p_{(\lambda)}(t)) = D_\infty(p_\infty(s'),p_\infty(t')).$$
\end{lemma}

\proof First notice that the condition $s,t\in [0,\gamma_\lambda(A)]\cup [\eta_\lambda(A),T]$
is equivalent to saying that $s',t'\in [-\eta_\infty(A),\gamma_\infty(A)]$
(recall that $\gamma_\infty(A)=\gamma_\lambda(A)$ and $\eta_\infty(A)= T-\eta_\lambda(A)$
on $\mathcal{F}_\lambda)$.

Then let $s,t\in [0,\gamma_\lambda(A)]\cup [\eta_\lambda(A),T]$ such that
$D^*_{\lambda}(\rho_{\lambda},p_{(\lambda)}(s))\leq r$ and $D^*_{\lambda}(\rho_{\lambda},p_{(\lambda)}(t))\leq r$.
By Lemma \ref{keyasymp}, $D^*_{\lambda}(p_{(\lambda)}(s),p_{(\lambda)}(t))$ is given by formula (\ref{tech10}). We 
claim that the right-hand side of this formula coincides with the 
right-hand side of formula (\ref{tech11}). To see this, let 
$s_0,t_0,s_1,t_1,\ldots,s_p,t_p\in [0,T]$ such that $s_0=s$
and $t_p=t$. For every $i\in\{0,1,\ldots,p\}$, set $s'_i=s_i$ if $s_i\leq T/2$
and $s'_i=s_i-T$ otherwise, and define $t'_i$ analogously. Then $s_0,s_1,\ldots,t_p\in [0,\gamma_\lambda(A^2)]\cup [\eta_\lambda(A^2),T]$
if and only if $s'_0,s'_1,\ldots,t'_p\in [-\eta_\infty(A^2),\gamma_\infty(A^2)]$. 
Assume that this condition holds. Then, one immediately checks that, for every $i\in\{0,1,\ldots,p\}$, 
we have $p_{(\lambda)}(s_i)=p_{(\lambda)}(t_i)$ if and only if $p_\infty(s'_i)=p_\infty(t'_i)$. 
Finally, we have also $ \wt D^\circ_{\lambda}(t_{i-1},s_i)= D^\circ_{\infty}(t'_{i-1},s'_i)$
for every $i\in\{0,1,\ldots,p\}$ (at this point it is crucial that $ \wt D^\circ_{\lambda}$
appears instead of $ D^\circ_{\lambda}$ in (\ref{tech10})). Our claim 
follows.

We cannot immediately infer that the right-hand side of (\ref{tech11})
coincides with $D_\infty(p_\infty(s),p_\infty(t))$ since we do not know yet that
$D_\infty(\rho_\infty,p_\infty(s))\leq r$ and $D_\infty(\rho_\infty,p_\infty(t))\leq r$.
However, the right-hand side of (\ref{tech11}) is clearly an upper bound for 
$D_\infty(p_\infty(s),p_\infty(t))$, and thus, by considering the special cases $s=0$ or $t=0$, we deduce 
from the equality between the right-hand sides
of (\ref{tech10}) and (\ref{tech11}) that $D_\infty(\rho_\infty,p_\infty(s))\leq r$
and $D_\infty(\rho_\infty,p_\infty(t))\leq r$. 

From a symmetric argument, we obtain that
the latter conditions imply $D^*_{\lambda}(\rho_{\lambda},p_{(\lambda)}(s))\leq r$ 
and $D^*_{\lambda}(\rho_{\lambda},p_{(\lambda)}(t))\leq r$. Finally, the equality 
between the right-hand sides
of (\ref{tech10}) and (\ref{tech11}) gives the identity
$D^*_{\lambda}(p_{(\lambda)}(s),p_{(\lambda)}(t)) = D_\infty(p_\infty(s'),p_\infty(t'))$. \endproof

\medskip
To complete the proof of Proposition \ref{couplingBMBP}, we verify that
the identity (\ref{idenball}) holds on $\mathcal{F}_\lambda$. Write $\bp_{(\lambda)}$ for
the composition of $p_{(\lambda)}$ with the canonical projection from
$\t_{(\lambda)}$ onto $Y^\lambda$, and similarly write 
$\bp_\infty$ for
the composition of $p_{\infty}$ with the canonical projection from
$\t_\infty$ onto $\pp$. Assuming that $\mathcal{F}_\lambda$ holds, we construct an isometry $\mathcal I$ from
$B_r(\mathbf{Y}^\lambda)$ onto $B_r(\boldsymbol{ \mathcal{P}})$, which maps $\rho_{\lambda}$ to $\rho_\infty$,
in the following way. An arbitrary point of $B_r(\mathbf{Y}^\lambda)$ is of the form 
$\bp_{(\lambda)}(s)$, where $s$ has to belong to
$[0,\gamma_\lambda(A)]\cup [\eta_\lambda(A),T]$ (by the first assertion 
in Lemma \ref{keyasymp} (i)). We then define
$${\mathcal I}(\bp_{(\lambda)}(s)) = \bp_\infty(s')$$
where $s'=s$ if $s\in [0,\gamma_\lambda(A)]$ and $s'=s-T$ if $s\in[\eta_\lambda(A),T]$,
as previously. The last assertion of Lemma \ref{corokey} shows that
this definition does not depend on the choice of $s$, and then that 
$\mathcal I$ is an isometry. Finally, using the first assertion of 
Lemma \ref{keyasymp} (ii), we easily get that $\mathcal I$
maps $B_r(\mathbf{Y}^\lambda)$ onto $B_r(\boldsymbol{ \mathcal{P}})$, and it is also immediate
that ${\mathcal I}(\rho_{\lambda})=\rho_\infty$. We conclude that the
pointed spaces $B_r(\mathbf{Y}^\lambda)$ and $B_r(\boldsymbol{ \mathcal{P}})$
are isometric on the event $\mathcal{F}_{\lambda}$. Since $\lambda\cdot \mathbf{m}_\infty$ has the same distribution
as $\mathbf{Y}^\lambda$, this completes the proof of Proposition \ref{couplingBMBP}
and of Theorem \ref{tangentcone}.

\section{The Brownian plane as the limit of discrete quadrangulations}

In this section, we prove Theorem \ref{scaling-limit}, which shows that the Brownian plane also arises as a scaling limit of discrete quadrangulations. Let us briefly discuss the strategy of the proof, which is, roughly speaking, the same as in the proof of Theorem \ref{tangentcone}. After recalling the constructions of discrete quadrangulations from labeled trees, we establish in Proposition \ref{comptree} a comparison lemma for these trees, which plays the same role as Proposition \ref{absocont} in the continuous setting. Using this lemma, we show in Proposition \ref{compmap} that we can couple the realizations of $Q_{n}$ and of $Q_{\infty}$ in such a way that the balls of radius $o(n^{1/4})$ centered at the origin are the same in $Q_{n}$ and $Q_{\infty}$ with large probability (this is the discrete counterpart to Proposition \ref{couplingBMBP}). We then use the convergence towards the Brownian map \eqref{cvbm} and Proposition \ref{couplingBMBP} to complete the proof. 

 \subsection{Discrete quadrangulations}
\label{disquad}

As in the introduction, we let $Q_{n}$ be uniformly distributed over the set of all rooted quadrangulations with $n$ faces.  Let us recall the definition of the uniform infinite planar quadrangulation (UIPQ).  More details can be found in \cite{CMM,Kri,Men}. 

If $m$ is a rooted planar map and $r \geq 0$, we define the ``combinatorial ball'' ${ \mathrm{Ball}}_{r}(m)$ as the submap of $m$ obtained by keeping only those 
edges and vertices of $m$ that are incident to (at least)  one face of $m$ having a vertex whose graph distance
from the root vertex is  smaller than or equal to $r$. Independently of the embedding chosen for $m$, this defines a planar map ${ \mathrm{Ball}}_{r}(m)$, which by convention has the same root as $m$.

Krikun \cite{Kri} proved that there exists a random infinite (rooted) planar quadrangulation $Q_{\infty}$ such that, for every $r \geq 0$, we have the following convergence in distribution
 \begin{eqnarray} \mathrm{Ball}_{r}(Q_{n}) &\xrightarrow[n\to \infty]{(d)}& \mathrm{Ball}_{r}(Q_\infty).\label{def:UIPQ}  \end{eqnarray} 
 We refer to \cite{CMM} for a discussion and a precise definition of infinite planar maps. The random infinite quadrangulation $Q_{\infty}$ is called the uniform infinite planar quadrangulation or UIPQ. Notice that a similar object has been defined earlier in the context of triangulations by Angel and Schramm, see \cite{Ang,AS}.

 An alternative approach to the preceding convergence 
 can be found in \cite{Men}, where it is also proved that $Q_\infty$ coincides with the random infinite quadrangulation that had been constructed 
 by Chassaing and Durhuus \cite{CD} from a random infinite well-labeled tree -- a different but related version of this construction, due to \cite{CMM},
 will be presented below.

As we already mentioned after the statement of Theorem \ref{scaling-limit}, it will be convenient to see discrete planar maps as length spaces: If $Q$ is a (finite or infinite) quadrangulation, we associate with $Q$ a pointed (boundedly compact) length space, denoted by $ \mathbf{Q}$, which is obtained by replacing every edge of $Q$ by a unit length Euclidean segment and pointing the resulting metric space at the root vertex of $Q$. More precisely, $ \mathbf{Q}$ is the union of 
 a (finite or infinite) collection of copies of the interval $[0,1]$ and two of these copies may intersect only at their endpoints if
 the associated edges of $Q$ share one or two vertices. Obviously, the distance between two points of $\mathbf{Q}$ is the length 
 of a shortest path between them. If $Q$ is finite, it is straightforward to verify that $ \op{d_{GH}}(V(Q), \mathbf Q)\leq 1$
 (here and later, we view $V(Q)$ as a pointed metric space, for the graph distance and with the root vertex of $Q$
 as distinguished point). More generally,
 both in the finite and the infinite case, we have $ \op{d_{GH}}(B_r(V(Q)), B_r(\mathbf Q))\leq 1$ for every $r\geq 0$. 
 
  We write $\mathbf{Q}_n$, resp.  $\mathbf{Q}_\infty$, for the pointed length space associated with $Q_n$,
  resp. with $Q_\infty$.

\subsection{Constructing quadrangulations from discrete trees}
 In this section we briefly recall the construction of discrete quadrangulations from labeled trees. The finite case is
presented in Section 5 of \cite{LGM}  and its extension to the infinite case can be found in \cite{CMM}
 (note that a different version of the construction in the infinite case had appeared earlier in \cite{CD}). The reader may consult the preceding references for more details.

\subsubsection{Labeled trees}

We start by recalling the standard formalism for plane trees. Let 
 \begin{eqnarray*}
\mathcal{U} &:=& \bigcup_{n=0}^{\infty} \mathbb{N}^n
,  \end{eqnarray*}
where $\mathbb{N} = \{ 1,2, \ldots \}$ and $\mathbb{N}^0 = \{
\varnothing \}$ by convention. An element $u$ of $\mathcal{U}$ is thus a 
(possibly empty) word made of positive integers, and $|u|\geq 0$ stands for
the length of $u$, sometimes also called the height of $u$.
If $u, v \in \mathcal{U}$,
$uv$ denotes the concatenation of $u$ and $v$. If $v\in \mathcal{U}\backslash\{\varnothing\}$ we can write 
$v=uj$, where $u\in\mathcal{U}$ and $j \in \mathbb{N}$, and we say that $u$ is the \emph{parent} of
$v$ or that $v$ is a \emph{child} of $u$. More generally, if $v$ is of
the form $uw$, for $u,w \in \mathcal{U}$, we say that $u$ is an
\emph{ancestor} of $v$ or that $v$ is a \emph{descendant} of $u$. A 
\emph{plane tree} $\tau$ is a (finite or infinite) subset of
$\mathcal{U}$ such that
\begin{enumerate}
\item $\varnothing \in \tau$ ($\varnothing$ is called the \emph{root}
  of $\tau$);
\item if $v \in \tau$ and $v \neq \varnothing$, the parent of $v$
  belongs to $\tau$;
\item for every $u \in \mathcal{U}$ there exists an integer $k_u(\tau) \geq 0$
  such that, for every $j\in\mathbb{N}$, $uj \in \tau$ if and only if $j \leq k_u(\tau)$.
\end{enumerate}

If the plane tree $\tau$ is finite, 
$|\tau|=\#(\tau)-1$ denotes the number of edges of $\tau$ and is called
the size of $\tau$.  

A \emph{ray} of an infinite plane tree $\tau$ is an infinite sequence
$u_{0},u_{1},u_{2},\ldots$ in $\tau$ such that $u_{0}=\varnothing$ and $u_i$
is the parent of $u_{i+1}$ for every $i\geq 0$. If an infinite tree $\tau$ has a unique ray, we call it the \emph{spine} of $\tau$ and denote it by $ \mathrm{S}_{\tau}(0), \mathrm{S}_{\tau}(1), \mathrm{S}_{\tau}(2), ...$. 

 
\medskip

In what follows, we say tree rather than plane tree. For every integer $n\geq 0$, we let $\mathbb{T}_n$
denote the set of all trees with size $n$. The cardinality of $\mathbb{T}_n$
is the Catalan number of order $n$,
$$\op{Cat}(n) = \frac{1}{n+1}{2n\choose n}.$$
For every $n\geq 0$, we let $T_{n}$ be uniformly distributed over $\mathbb{T}_n$.

\begin{figure}[!h]
 \begin{center}
 \includegraphics[width=8cm]{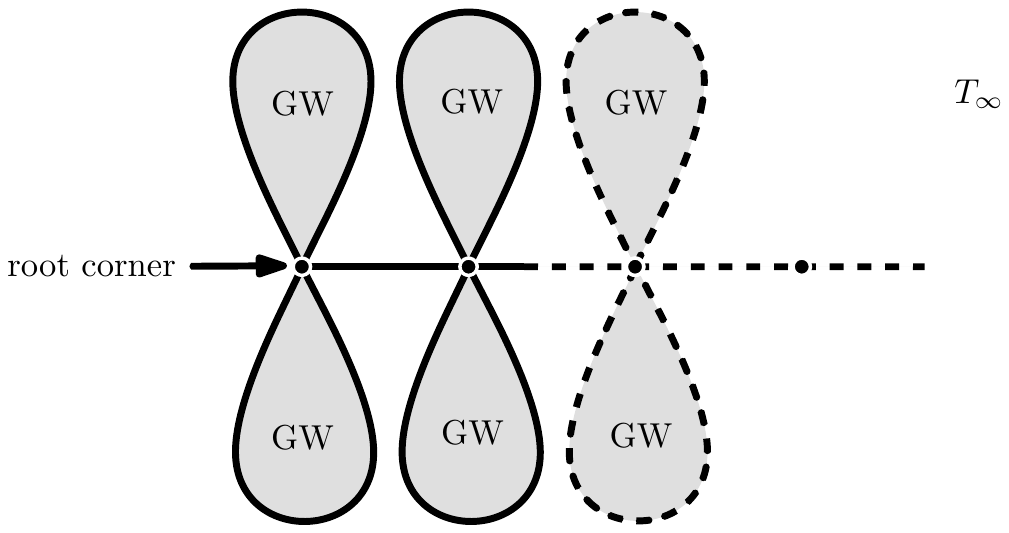}
 \caption{ \label{tinfty}Construction of $T_{\infty}$. The  trees that are grafted on the spine are
 independent Galton-Watson trees with geometric offspring distribution with parameter $1/2$.}
 \end{center}
 \end{figure}

\medskip 
\paragraph{The uniform infinite tree $T_{\infty}$.} We now introduce the infinite ``local limit'' of the 
random trees $T_{n}$ as $n \to \infty$. If $\tau$ is a plane tree and $k \in \{0,1,2, \ldots \}$, we let $[\tau]_{k}:=\{ v \in \tau : |v| \leq k \}$ be the plane tree obtained from $\tau$ by keeping only its first $k$ generations. It follows from the work of Kesten \cite{Kes} (see also \cite{LPP}) that there exists an infinite random tree $T_{\infty}$ with only one spine such that, for every $k \geq 0$, we have 
 \begin{eqnarray} \label{convergence} [T_{n}]_{k} & \xrightarrow[n\to\infty]{(d)} & [T_{\infty}]_{k}.  \end{eqnarray}
The tree $T_{\infty}$ is known as the critical geometric Galton-Watson tree conditioned to survive. It can be described informally as follows (a rigorous 
presentation using the preceding formalism for trees can be found in \cite{CMM}): Start with a semi-infinite line of vertices that will be the spine of $T_{\infty}$ and graft to the left and to the right of each vertex of the spine independent  critical geometric Galton-Watson trees (with parameter $1/2$). The root of $T_\infty$ is obviously the first vertex of the spine. See Fig.\,\ref{tinfty}.

\paragraph{Labeled trees.} A \emph{ labeled tree} is a pair $\theta = (\tau, (\ell(u))_{u \in \tau})$ that consists of a  tree $\tau$ and a collection of integer labels assigned to the vertices of $\tau$, such that $\ell(\varnothing)=0$ and $|\ell(u) - \ell(v)| \leq 1$ whenever $u,v \in \tau$ are neighbors
(meaning that $u$ is either the parent or a child of $v$). We denote the size of $\theta$ by $|\theta| = |\tau|$.  We let $ \mathcal{T}$ stand for the set of all finite labeled trees and we let  $ \mathscr{S}$ be the set of all infinite labeled trees $\theta=(\tau,\ell)$ such that $\tau$ has only one spine and $\inf_{i \geq 0} \ell(\mathrm{S}_{\tau}(i)) = -\infty$.\medskip

If $\tau$ is a (finite or infinite)  tree, we can assign labels to its vertices in a uniform way: Write $E_\tau$ for the set of all edges of $\tau$ and consider a collection $(\eta_e)_{e\in E_\tau}$ of independent random variables uniformly distributed over $\{-1,0,+1\}$. For any vertex $u \in \tau$, the 
(random) label $\ell(u)$ of $u$ is then defined as the sum  of the variables $\eta_e$ over all edges $e\in E_\tau$ belonging to the geodesic path from the root to $u$. In particular $\ell(\varnothing)=0$. The resulting random labeled tree $\theta=(\tau,\ell)$  is called the uniform labeling of $\tau$. 
If $\tau$ is random, we can still consider its uniform labeling by applying the preceding construction after conditioning on $\tau$. This applies in particular to the random tree  $T_{n}$, and the resulting random labeled tree $\Theta_{n}= (T_{n}, \ell_{n})$ is uniformly distributed over the set
$$\mathcal{T}_n:=\{ \theta \in \mathcal{T} : | \theta|=n\}.$$ 
If we apply the construction to  $T_{\infty}$, the resulting labeled tree $\Theta_{\infty}=(T_{\infty}, \ell_{\infty})$ is called the uniform infinite labeled tree.
Notice that $\Theta_{\infty} \in \mathscr{S}$ almost surely.

\subsubsection{The construction in the finite case} We will describe Schaeffer's bijection between $\mathcal{T}_n \times\{0,1\}$ and the set of all rooted and pointed quadrangulations with $n$ faces. See Section 5 of
\cite{LGM} for more details.

\begin{figure}[!h]
 \begin{center}
 \includegraphics[width=6cm]{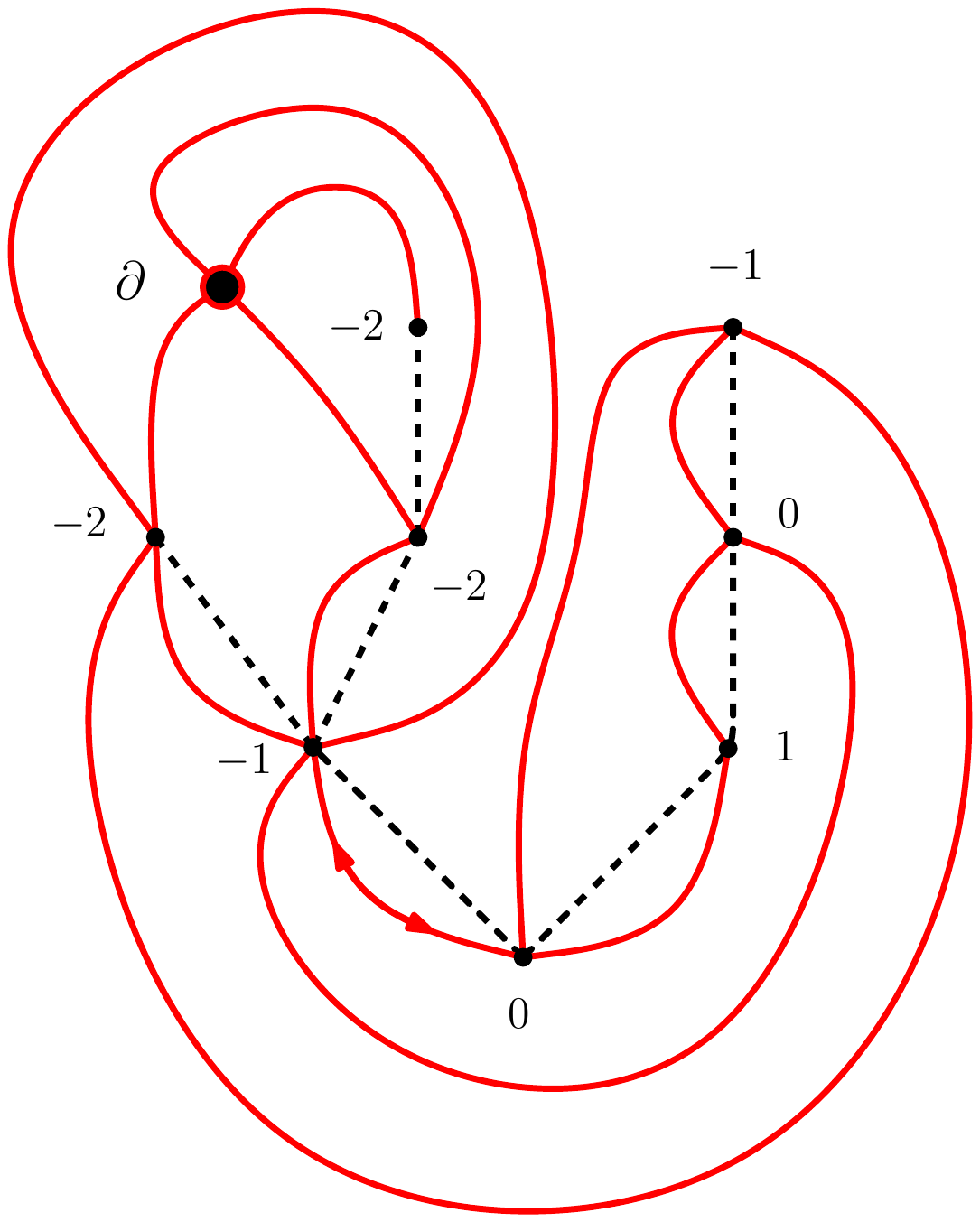}
 \caption{ \label{schaeffer-fini} Construction of a rooted and pointed quadrangulation from a
 labeled tree.}
 \end{center}
 \end{figure}

A rooted and pointed quadrangulation is a pair $(q,\partial)$ consisting of a rooted quadrangulation $q$ together with a distinguished vertex $\partial$ of $q$.
We start from a finite labeled tree $\theta=(\tau,\ell)\in\mathcal{T}_n$ and write $\min\ell$ for the minimal label on the tree. In order to associate a rooted and pointed quadrangulation with $\theta$, we first embed the tree $\tau$  in the plane, in the way suggested by Fig.\ref{schaeffer-fini} (on this figure the edges of $\tau$ are the dotted lines, the figures are the labels of vertices, the root vertex
is obviously at the bottom and the lexicographical order between children of a given vertex corresponds to listing the edges 
from the left to the right). We also add an extra vertex (outside the embedded tree) denoted by $\partial$. A corner $c$ of the (embedded) tree is  an angular sector between two adjacent edges. The label $\ell(c)$ of the corner $c$ is the label of the associated vertex $ \mathcal{V}(c)$. Corners of the tree have a cyclic ordering given by the clockwise contour of the tree in the plane (if we imagine a particle that moves around the embedded tree along its edges
in clockwise order,
it will visit every corner exactly once before coming back to the corner it started from). 

To obtain the edges of the quadrangulation $q$ associated with $\theta$, we proceed as follows. For every corner $c$ of the tree
such that $\ell(c)>\min\ell$, we draw an edge between $c$ and the first corner after $c$ with label $\ell(c)-1$. This corner is called the successor of $c$ and denoted by $ \mathcal{S}(c)$. If $\ell(c)=\min\ell$ (so that the definition of $\mathcal{S}(c)$
does not make sense),  we draw an edge between $c$ and $\partial$. This construction can be made in such a way that the edges 
do not cross each other and do not cross the edges of the tree $\tau$ (see Fig.\ref{schaeffer-fini}). After erasing the embedding of the tree, the resulting map $q$ is a quadrangulation with $n$ faces and whose vertices are
exactly the vertices of $\tau$ plus the extra vertex $\partial$. Futhermore,  the labeling has the following interpretation: For every $u \in \tau$, we have 
 \begin{eqnarray} \mathrm{d_{gr}}( \partial, u) &=& \ell(u) - \min \ell +1,  \label{distances}\end{eqnarray} where $ \mathrm{d_{gr}}$ stands for the graph distance in the quadrangulation. 
 
The quadrangulation $q$ is pointed at the vertex $\partial$ and its root edge is the edge of $q$ which is drawn from the root corner $c_0$ of
$\tau$ (the root corner is the corner ``below'' the root $\varnothing$ of $\tau$).  In order to specify the orientation of the root edge, we need an extra variable $\epsilon \in \{0,1\}$: The origin of the root edge is the vertex $\varnothing$ if $\epsilon =0$ and the other end of the root edge 
if $\epsilon=1$. 

Call $\Phi(\theta,\epsilon)$ the rooted and pointed quadrangulation that is obtained by the preceding construction. Then $\Phi$
is a bijection between $\mathcal{T}_n \times\{0,1\}$ and the set of all rooted and pointed quadrangulations with $n$ faces. 
Consequently,  if $\Theta_n$ is a uniform labeled tree with $n$ edges and $\eta$ is an independent Bernoulli variable of parameter $1/2$, then 
the random rooted quadrangulation derived from $\Phi(\Theta_n,\eta)$ by ``forgetting'' the distinguished vertex has the same 
distribution as $Q_n$. 
  
\subsubsection{The construction in the infinite case}  

\begin{figure}[!h]
 \begin{center}
 \includegraphics[width=10cm]{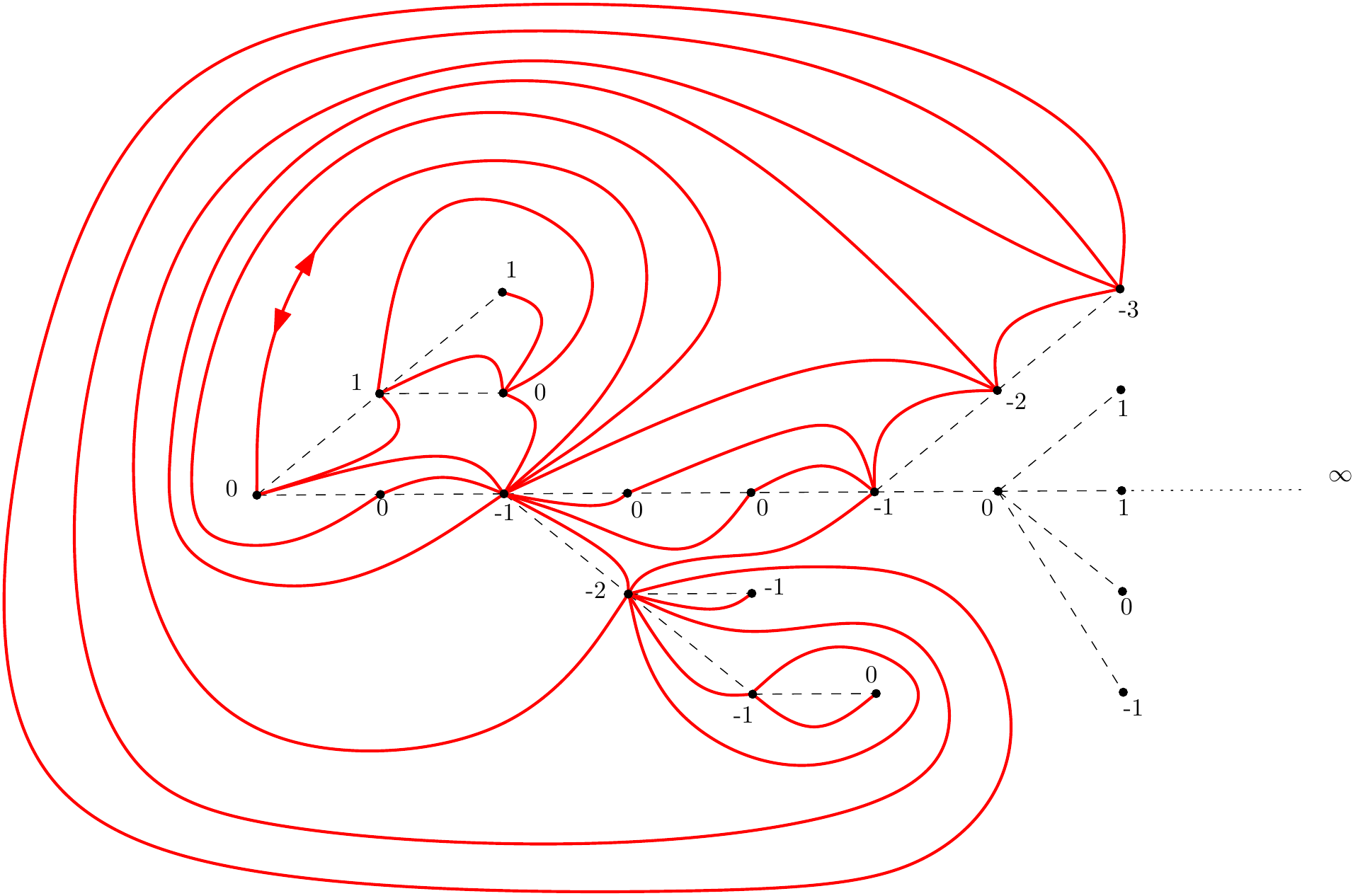}
 \caption{ \label{Schaeffer-infini} Construction of a rooted infinite quadrangulation from a
tree of $\mathscr{S}$.}
 \end{center}
 \end{figure}

The preceding construction can easily be extended to the case when the tree $\theta=(\tau,\ell)$ is an infinite labeled tree in $ \mathscr{S}$. We first embed the infinite tree properly in the plane and then draw an edge between each corner $c$ of the tree and the first corner  with label $\ell(c)-1$ following $c$ in the clockwise contour order. We again denote this corner by $ \mathcal{S}(c)$ and call it the successor of $c$. Because of our assumption $\theta\in \mathscr{S}$, there is no vertex of minimal label and thus unlike the finite case there is no need to add an extra vertex $\partial$. The construction should be clear from Fig.\ref{Schaeffer-infini}. The resulting planar map in an infinite quadrangulation, see \cite{CMM}. We root this quadrangulation at the edge connecting  the root corner of $\tau$ to its successor. As in the finite case, the orientation of this edge is given by an extra variable $\epsilon \in \{0,1\}$. This infinite rooted quadrangulation is denoted by $\Phi'(\theta,\epsilon)$. Note that the vertex set of $\Phi'(\theta,\epsilon)$ is precisely the vertex set of $\tau$. 

Finally, from \cite[Theorem 1]{CMM} we know that if $\Theta_{\infty}$ is a uniform infinite labeled tree and $\eta$ is an independent Bernoulli variable of parameter $1/2$ then 
    \begin{eqnarray*} \Phi'(\Theta_{\infty},\eta) &\overset{(d)}{=}& Q_{\infty}. \end{eqnarray*}

\subsection{Comparison lemmas}

\subsubsection{Comparison of trees}
One of the key ingredients in  the proof of Theorem \ref{scaling-limit} is an improvement of the convergence \eqref{convergence}. Roughly speaking, we will see that the fixed integer $k \geq 0$ in \eqref{convergence} can be replaced by a sequence $k=k(n)$ as soon as $k(n) = o( \sqrt{n})$ as $n \to \infty$. For technical reasons, it is more convenient to deal with \emph{pointed} trees.  \bigskip

\begin{figure}[!h]
\begin{center} 
\includegraphics[width=12cm]{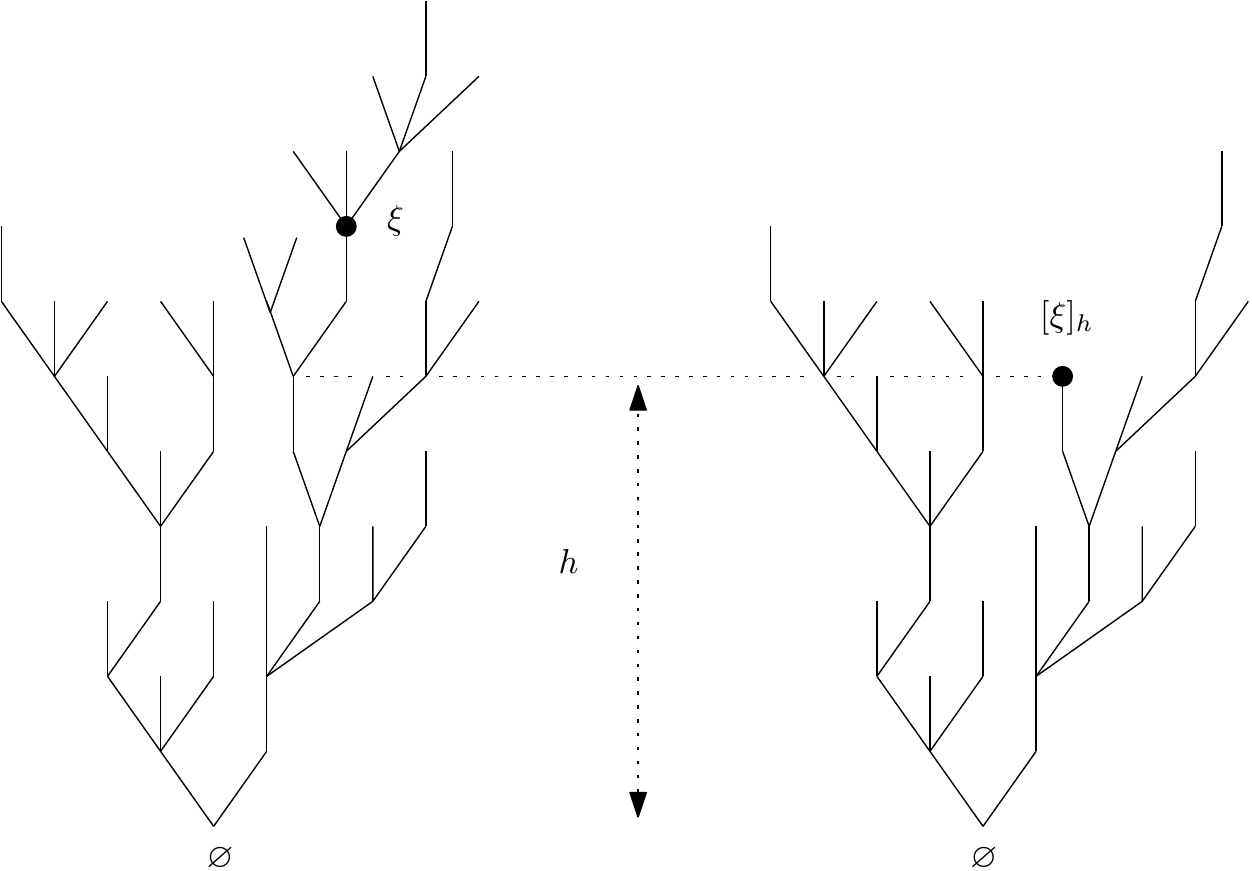}
\caption{ \label{fig:pruned}A pointed tree  $\boldsymbol{\tau}=(\tau,\xi)$ and the resulting pruned tree $ \mathscr{P}( \boldsymbol{\tau},h)$ at height $h$.}\end{center}
\end{figure}

A (finite) pointed tree is a pair $ \boldsymbol{\tau}=(\tau,\xi)$ where $\tau$ is a finite
 tree and $\xi$ is a distinguished vertex of $\tau$. Let $ \boldsymbol{\tau}=(\tau,\xi)$ 
 be a pointed tree. For every integer $h$ such that $0 \leq h <
|\xi|$, we let $  \mathscr{P}( \boldsymbol{\tau},h)$ stand for the subtree of $\tau$ consisting of  all  vertices $u \in \tau$ such that the 
height of the most recent common ancestor of $u$ and $\xi$ is
strictly less than $h$, together with the ancestor of $\xi$ at
height exactly $h$, which is denoted by $[\xi]_{h}$. We furthermore point this tree at $[\xi]_{h}$. See Fig.\ref{fig:pruned}.
By convention when  $h \geq |\xi|$, we declare that
$\mathscr{P}(\boldsymbol{\tau},h) = (\{\varnothing\},\varnothing)$. 

If $\tau
$ is an infinite tree with only one spine (denoted as above by  $\mathrm{S}_{\tau}(0), \mathrm{S}_{\tau}(1), \ldots$), we can extend the former definition by considering $\tau$ as pointed at infinity.  
Formally, for every integer $h \geq 0$, we let $\mathscr{P}(\tau,h)$ be the subtree of $\tau$ consisting of the vertices $\mathrm{S}_{\tau}(0), \ldots , \mathrm{S}_{\tau}(h)$ 
of the spine, together with the subtrees grafted to the left and the right of $\mathrm{S}_\tau(i)$ for $0\leq i\leq h-1$. We point this tree at $ \mathrm{S}_{\tau}(h)$. 

If $ \boldsymbol{\tau}$ is a pointed tree, and $h,h'$ are integers such that $0\leq h\leq h'$, then one immediately checks that
\begin{equation}
\label{compopruning}
\mathscr{P}(   \mathscr{P}( \boldsymbol{\tau},h'),h)=   \mathscr{P}( \boldsymbol{\tau},h).
\end{equation}
The same property holds if we replace $\boldsymbol{\tau}$ by an infinite tree having only one spine. 
\medskip

In what follows,
$\mathbf{T}_{n}=(T_{n},\xi_{n})$ is uniformly distributed over the set of all pointed 
trees with size $n$. This is consistent with our previous notation, since $T_{n}$ is then uniformly distributed over $\mathbb{T}_n$. 
Note that, conditionally on $T_{n}$, the vertex $\xi_{n}$ is uniformly distributed over the vertex set of $T_n$. In particular if $ \boldsymbol{\tau_{0}}$ is a fixed pointed tree
with size $n$, we have
$$ {P}( \mathbf{T}_{n} = \boldsymbol{\tau_{0}}) =\frac{1}{(n+1)\op{Cat}(n)}
=\frac {1}{{2n \choose n}}.$$ 
Recall that $T_{\infty}$ denotes the uniform infinite tree. The following
proposition relates $ \mathbf{T}_{n}$ to $T_{\infty}$. 

\begin{proposition} \label{comptree}
For every $\varepsilon>0$, there exists $\delta>0$ 
such that, for every sufficiently large integer $n$, the bound
 \begin{eqnarray*} \bigg| {P}\Big(\mathscr{P}\big(\mathbf{T}_{n}, 
\lfloor\delta n^{1/2}\rfloor \big) \in A\Big) - {P} \Big(\mathscr{P}\big(
T_{\infty}, \lfloor\delta n^{1/2}\rfloor
\big) \in A \Big) \bigg| &\leq& \varepsilon, \end{eqnarray*} 
holds for any set $A$
of pointed trees. Consequently, if $k_{n}$ is a sequence of positive integers such that $k_n=o(n^{1/2})$
as $n\to\infty$, then the total variation distance between the law of $\mathscr{P}( \mathbf{T}_{n}, k_{n})$ and 
that of $\mathscr{P}(T_{\infty},k_{n})$ tends to $0$ as $n $ goes to $\infty$. \end{proposition}

\noindent{\bf Remark.} It is easy to see that Proposition \ref{comptree} does imply  \eqref{convergence} and is in fact much stronger. A similar result 
has been proved by Aldous \cite[Theorem 2]{Al1} for Poisson Galton-Watson trees.

\proof
Let $h \geq 0$ be an integer. We first
identify the distribution of $\mathscr{P}(T_{\infty},
h)$. Recall that  $\mathscr{P}(T_{\infty},h)$ consists of the
fragment of the spine of $T_{\infty}$ up to height $h$, together with
the subtrees grafted to the left and to the right of the spine up to height
$h-1$, and that this tree is pointed at $\mathrm{S}_{T_{\infty}}(h)$.  The subtrees grafted on the spine of $T_{\infty}$ are
independent Galton-Watson trees with geometric offspring distribution of parameter
$1/2$. If $\mathbb{Q}_{GW}$ stands for the distribution of one of these trees, a standard
formula for Galton-Watson trees gives, for every finite tree $\tau$,
\begin{equation}
\label{lawGWT}
\mathbb{Q}_{GW}(\tau) = \prod_{u\in\tau} 2^{-k_u(\tau)-1} = 2^{-2|\tau|-1}= \frac{1}{2} 4^{-|\tau|}.
\end{equation}
Let $ \boldsymbol{\tau_{0}}=(\tau_{0},\xi)$ be a pointed tree, such that $\xi$ is a vertex of $\tau_{0}$ at height $h$ with no child. Write $\tau_{(i)}$, $0\leq i\leq h-1$, resp. $\tau'_{(i)}$, $0\leq i\leq h-1$, for the successive subtrees that branch from the left side, resp. from
the right side, of the
ancestral line of $\xi$ in $\tau_0$. It follows from the previous observations that
\begin{equation}
\label{loiinfty} {P} \big( \mathscr{P}\big(T_{\infty} ,h \big)= \boldsymbol{\tau_{0}}\big) =
\prod_{i=0}^{h-1} \big(\frac{1}{2} 4^{-|\tau_{(i)}|}\big)\big(\frac{1}{2} 4^{-|\tau'_{(i)}|}\big)= 4^{-|\tau_{0}|}.
\end{equation}
We also notice that the size  of the random pointed tree $\mathscr{P}(T_{\infty},h)$ can be written in the form 
\begin{eqnarray*} |\mathscr{P}(T_{\infty},h)| &=& h + \sum_{i=0}^{h-1}  (\mathsf{N}_{i,\ell} +\mathsf{N}_{i,r}),   \end{eqnarray*} 
where the random variables 
$\mathsf{N}_{i,\ell}$, $\mathsf{N}_{i,r}$, for $0\leq i\leq h-1$ are independent and distributed according to the
size of a Galton-Watson tree with geometric offspring distribution of parameter
$1/2$.  For every integer $n\geq0$, \eqref{lawGWT} gives
\begin{eqnarray*}
\label{size} {P}\left(\mathsf{N}_{0, \ell} = n \right)  = \frac{1}{2}\op{Cat}(n) 4^{-n} \quad \underset{n \to \infty}{\sim} \quad \frac{n^{-3/2}}{2\sqrt{\pi}}.\end{eqnarray*}  Standard facts about domains of attraction (see for example \cite[p. 343-350]{BGT}) thus imply that $n^{-2} |\mathscr{P}(T_{\infty},n)|$ converges in distribution towards a stable law with parameter $1/2$ (we could also derive this from the fact that
$2\mathsf{N}_{0, \ell} +1$ is distributed as the hitting time of $1$ for simple random
walk on $\mathbb{Z}$ started from $0$). In particular, if $\varepsilon \in (0,1)$ is fixed, there exists 
a constant $c_{\varepsilon}>1$ such that for every integer $k \geq 0$ we have
\begin{eqnarray}
\label{tailleinfty} {P}\Big( \big|\mathscr{P}(T_{\infty},k)\big| \leq c_{\varepsilon}k^2 \Big) & \geq & 1-\varepsilon.\end{eqnarray}

We now compute the distribution of
$\mathscr{P}(\mathbf{T}_{n},h)$. Let  the pointed tree $ \boldsymbol{\tau_{0}}=(\tau_{0},\xi)$ be as previously, with $|\xi|=h$. The event
$\{\mathscr{P}( \mathbf{T}_{n},h)= \boldsymbol{\tau_{0}}\}$ holds if and only
if the tree $T_{n}$ is obtained from the tree $\tau_{0}$ by grafting at
$\xi$ a subtree $ \mathfrak{t}$ having $n-| \tau_{0}|$ edges and if furthermore
the distinguished point $\xi_{n}$ of $T_{n}$ is in $ \mathfrak{t}$
but  is different from its root. Hence a direct counting  argument shows that
$$ {P}\big(\mathscr{P}(\mathbf{T}_{n},h)= \boldsymbol{ \tau}_{0}\big)  \quad=\quad 
\mathbf{1}_{\{|\tau_0|<n\}}\; \frac{(n-| \tau_{0}|)\op{Cat}(n-| \tau_{0}|)}{(n+1)\op{Cat}(n)}.$$
Recall that we have fixed $\varepsilon\in(0,1)$. Using asymptotics for Catalan numbers, we can find an integer $N_\varepsilon\geq 1$, which does not
depend on $h$ nor on the choice of the pointed tree $ \boldsymbol{\tau_{0}}=(\tau_{0},\xi)$ satisfying the preceding properties, 
such that, for every integer $n\geq |\tau_0| + N_\varepsilon$,
$$ (1-\varepsilon)\,4^{-| \tau_{0}|}\left( 1- \frac{|  \tau_{0}|}{n}\right)^{-1/2}\leq
 {P}\big(\mathscr{P}(\mathbf{T}_{n},h)= \boldsymbol{ \tau}_{0}\big) 
 \leq (1+\varepsilon)\,4^{-| \tau_{0}|}\left( 1- \frac{|  \tau_{0}|}{n}\right)^{-1/2},
$$
 and therefore, using (\ref{loiinfty}),
\begin{equation}
 \label{equivalent}  
(1-\varepsilon)\left( 1- \frac{|  \tau_{0}|}{n}\right)^{-1/2}\leq
 \frac{{P}\big(\mathscr{P}(\mathbf{T}_{n},h)= \boldsymbol{ \tau}_{0}\big) }{{P} \big( \mathscr{P}\big(T_{\infty} ,h \big)= \boldsymbol{\tau_{0}}\big)}
 \leq (1+\varepsilon)\left( 1- \frac{|  \tau_{0}|}{n}\right)^{-1/2}.
 \end{equation}
Recall the constant $c_\varepsilon$ from \eqref{tailleinfty}. We choose $\delta >0$ small enough so that $c_\varepsilon\delta^2<1/2$
and $(1-c_\varepsilon\delta^2)^{-1/2}<1+\varepsilon$. We apply (\ref{equivalent}) with $h=\lfloor \delta n^{1/2}\rfloor$, and we get
that, for every $n\geq 2 N_\varepsilon$, and for every pointed tree $ \boldsymbol{\tau_{0}}=(\tau_{0},\xi)$ 
with $|\tau_0|\leq c_\varepsilon\delta^2n$, such 
that $|\xi|= \lfloor \delta n^{1/2}\rfloor$ and $\xi$ has no child, 
$$ 1-\varepsilon\leq 
 \frac{{P}\big(\mathscr{P}(\mathbf{T}_{n},\lfloor \delta n^{1/2}\rfloor)= \boldsymbol{ \tau}_{0}\big) }
 {{P} \big( \mathscr{P}\big(T_{\infty} ,\lfloor \delta n^{1/2}\rfloor \big)= \boldsymbol{\tau_{0}}\big)}\leq (1+\varepsilon)^2.$$
 By (\ref{tailleinfty}), we have ${P}\big( \big|\mathscr{P}(T_{\infty},\lfloor \delta n^{1/2}\rfloor)\big| \leq c_{\varepsilon}\delta^2n \big) \geq 1-\varepsilon$, and 
 it then follows from the preceding bounds that, for every $n\geq 2N_\varepsilon$,
 $${P}\Big( \big|\mathscr{P}(\mathbf{T}_{n},\lfloor \delta n^{1/2}\rfloor)\big| \leq c_{\varepsilon}\delta^2n \Big) \geq (1-\varepsilon)^2.$$
 Finally, if $n\geq 2N_\varepsilon$ and $A$ is any set of pointed trees,
  \begin{eqnarray*} &&\bigg| {P}\Big(\mathscr{P}\big(\mathbf{T}_{n}, 
\lfloor\delta n^{1/2}\rfloor \big) \in A\Big) - {P} \Big(\mathscr{P}\big(
T_{\infty}, \lfloor\delta n^{1/2}\rfloor
\big) \in A \Big) \bigg| \\
&&\quad\leq{P}\Big( \big|\mathscr{P}(T_{\infty},\lfloor \delta n^{1/2}\rfloor)\big| > c_{\varepsilon}\delta^2n \Big)
+ {P}\Big( \big|\mathscr{P}(\mathbf{T}_{n},\lfloor \delta n^{1/2}\rfloor)\big| >  c_{\varepsilon}\delta^2n \Big) \\
&&\qquad + \sum_{ \boldsymbol{\tau_{0}}\in A\,,\, |{\tau_{0}}|\leq  c_{\varepsilon}\delta^2n} 
 {P} \big( \mathscr{P}\big(T_{\infty} ,\lfloor \delta n^{1/2}\rfloor \big)= \boldsymbol{\tau_{0}}\big)
 \left|  \frac{{P}\big(\mathscr{P}(\mathbf{T}_{n},\lfloor \delta n^{1/2}\rfloor)= \boldsymbol{ \tau}_{0}\big) }
 {{P} \big( \mathscr{P}\big(T_{\infty} ,\lfloor \delta n^{1/2}\rfloor \big)= \boldsymbol{\tau_{0}}\big)} -1 \right|\\
 &&\quad\leq \varepsilon + (1-(1-\varepsilon)^2) + (1+\varepsilon)^2 -1\\
&&\quad\leq 5\varepsilon. \end{eqnarray*} 
 This completes the proof of the first assertion of the proposition.
 The second assertion easily follows from the first one by using the composition property (\ref{compopruning}). 
\endproof

\subsubsection{Comparison of maps}
In this section we use Proposition \ref{comptree} to derive a version of the convergence \eqref{def:UIPQ}
where the radius $r$ can tend to infinity with $n$. Recall that if $Q$ is a (finite or infinite) rooted quadrangulation, 
the vertex set $V(Q)$ of $Q$ is viewed as a pointed metric space. To simplify notation,  we write $B_r(Q)=B_r(V(Q))$. 
Note that $B_r(Q)$ is again a metric space pointed at the root vertex of $Q$.

\begin{lemma} \label{geometric} Let $ \theta=(\tau,\ell)$ be a finite labeled tree and let $\theta'=( \tau',\ell')$ be an infinite labeled tree in $ \mathscr{S}$. 
Let $\eta\in\{0,1\}$ and 
let $Q=\Phi(\theta,\eta)$, resp. $Q'=\Phi'(\theta',\eta)$, be the rooted quadrangulation, resp. the rooted infinite quadrangulation, constructed
from the pair $(\theta,\eta)$, resp. from the pair $(\theta',\eta)$, via Schaeffer's bijection. 
Assume that there exist $\xi \in \tau$ and an integer $h \geq 1$ such that $ \mathscr{P}((\tau,\xi),h) = \mathscr{P}(\tau',h)$ and $\ell(u)=\ell'(u)$ 
for every $u\in \mathscr{P}(({\tau},\xi),h)$. Set 
$$r \ =\ -\min_{0 \leq i \leq h} \ell\big(\mathrm{S}_{\tau'}(i)\big)\geq 0\,.$$
Assume that $r\geq 3$ and set $r'=\frac{1}{2}(r-3)$.
Then,
$$ B_{r'}\big({Q}\big) = B_{r'}\big({Q'}\big).$$
\end{lemma}
 
\noindent{\bf Remark.} The conclusion of the lemma should be understood in the
 sense that $ B_{r'}\big({Q}\big)$ and $B_{r'}\big({Q'}\big)$ 
 are isometric as pointed metric spaces.

 \proof  We use the notation $\op{d}^Q_{\rm gr}$, resp. $\op{d}^{Q'}_{\rm gr}$ for the
 graph distance on $V(Q)$, resp. on $V(Q')$. Recall that $V(Q)$ can be identified with 
 the tree $\tau$ (plus an extra vertex that plays no role in this proof)
 and similarly $V(Q')$ is identified with $\tau'$.
 If $u,v\in \tau$, write $\llbracket u,v\rrbracket$ for the geodesic path between $u$ and
 $v$ in the tree $\tau$. A simple consequence of the construction of edges in Schaeffer's bijection
 is the discrete cactus bound (see formula (4) in \cite{CMM} and compare with \eqref{cabo} above) stating that
 $$\op{d}^Q_{\rm gr}(u,v) \geq \ell(u) + \ell(v) - 2 \min_{w\in \llbracket u,v \rrbracket} \ell(w).$$
 The same bound holds for $\op{d}^{Q'}_{\rm gr}(u,v)$ when $u,v\in \tau'$, replacing $\ell$ by $\ell'$. 
 
 Let $k\in\{0,1,\ldots,h\}$ be such that
 $$\ell(\mathrm{S}_{\tau'}(k))=\min_{0 \leq i \leq h} \ell\big(\mathrm{S}_{\tau'}(i)\big).$$
 Note that we have also $ \mathscr{P}((\tau,\xi),k) = \mathscr{P}(\tau',k)$ by our assumption and (\ref{compopruning}).
 We then observe that, if $u\in \tau\backslash \mathscr{P}(\tau,k)$, the ancestral line of 
 $u$ coincides with the ancestral line of $\xi$ at least up to level $k$,
 and thus must contain the vertex $\mathrm{S}_{\tau'}(k)$ (which belongs to the latter ancestral line
 by our assumption $ \mathscr{P}((\tau,\xi),h) = \mathscr{P}(\tau',h)$). The cactus bound then gives
 $$\op{d}^Q_{\rm gr}(\varnothing, u) \geq -\min_{w\in\llbracket \varnothing, u\rrbracket} \ell(w) \geq -\ell(\mathrm{S}_{\tau'}(k)) =r.$$
 If $u\in \tau'\backslash  \mathscr{P}(\tau',k)$, the same argument gives 
 $$\op{d}^{Q'}_{\rm gr}(\varnothing, u)\geq r.$$
 
 Now let $u\in   \mathscr{P}((\tau,\xi),k) = \mathscr{P}(\tau',k)$ such that $\op{d}^Q_{\rm gr}(\varnothing,u)\leq r-1$. Then any vertex $v$
 that belongs to a geodesic path from $\varnothing$ to $u$ in $Q$ must satisfy $\op{d}^Q_{\rm gr}(\varnothing,v)\leq r-1$
 and therefore also belong to $ \mathscr{P}((\tau,\xi),k)$. 
 We claim that any edge of $Q$ that appears on this geodesic path must correspond to an edge of $Q'$ with the same
 endpoints. This follows from the construction of edges in Schaeffer's bijection (both in the finite and in the infinite case),
 except that we must rule out the possibility of an edge starting from the left side of the ancestral line of $[\xi]_k=\mathrm{S}_{\tau'}(k)$
 in $ \mathscr{P}((\tau,\xi),k)$ and ending on the right side of this ancestral line (obviously such edges do not appear
 in $Q'$).  However, such an edge would start from a corner $c$ and end at a corner $c'$ such that the set of 
 all corners between $c$ and $c'$ (in cyclic ordering) would contain a corner of the vertex $\mathrm{S}_{\tau'}(k)$. 
 But then the label of both endpoints of the edge would be smaller than $\ell(\mathrm{S}_{\tau'}(k))= -r$. This is absurd 
 since labels on the geodesic path from $\varnothing$ to $u$ must be greater than or equal to $-r+1$ since
 $\op{d}^Q_{\rm gr}(\varnothing,u)\leq r-1$.
 
 The preceding discussion entails that if $u\in   \mathscr{P}((\tau,\xi),k) = \mathscr{P}(\tau',k)$ is such that $\op{d}^Q_{\rm gr}(\varnothing,u)\leq r-1$
 then  $\op{d}^{Q'}_{\rm gr}(\varnothing,u)\leq \op{d}^{Q}_{\rm gr}(\varnothing,u)$. But the same argument (in fact easier since we do not have to rule
 out the possibility of edges going from the left side of the spine to its right side) also shows 
 that  if $u\in   \mathscr{P}(\tau',k)$ is such that $\op{d}^{Q'}_{\rm gr}(\varnothing,u)\leq r-1$
 then  $\op{d}^{Q}_{\rm gr}(\varnothing,u)\leq \op{d}^{Q'}_{\rm gr}(\varnothing,u)$. Hence vertices that are at distance less than 
 $r-1$ from $\varnothing$ are the same in $Q$ and in $Q'$. 
 
We next observe that, if $u$ and $v$ are two vertices of $\mathscr{P}((\tau,\xi),k)$ such that $\op{d}^Q_{\rm gr}(\varnothing,u)\leq \frac{1}{2}(r-1)$
and $\op{d}^Q_{\rm gr}(\varnothing,v)\leq \frac{1}{2}(r-1)$, we have 
\begin{equation} 
\label{keygeometric}
\op{d}^Q_{\rm gr}(u,v)= \op{d}^{Q'}_{\rm gr}(u,v). 
\end{equation}
Indeed, any vertex $w$ on a geodesic path from $u$ to $v$ in $Q$ must be at $\op{d}^Q_{\rm gr}$-distance at most 
$\frac{1}{2}(r-1)$ from either $u$ or $v$, and thus at $\op{d}^Q_{\rm gr}$-distance at most $r-1$ from $\varnothing$.
By the same argument as previously, any edge on a geodesic path from $u$ to $v$ in $Q$
corresponds to an edge in $Q'$, and the converse also holds. The equality \eqref{keygeometric} follows.

Finally \eqref{keygeometric} and the preceding considerations show that the ball of radius $\frac{1}{2}(r-1)$
centered at $\varnothing$ in $Q$ is isometric to the same ball in $Q'$. Since the root vertex of both $Q$ and $Q'$ is
either $\varnothing$ (if $\eta=0$) or the successor of the first corner of $\varnothing$, which is at graph distance $1$
from $\varnothing$, the conclusion of the lemma follows.
 \endproof

Recall that $Q_{n}$ is uniformly distributed over the set of all rooted quadrangulations with $n$ edges and that $Q_{\infty}$ is  the uniform infinite planar quadrangulation.

\begin{proposition} \label{compmap} For every $\varepsilon>0$, there exists $\alpha>0$ 
such that, for every sufficiently large integer $n$ and every $m\in\{n+1,n+2,\ldots\}$, we can construct $Q_{n}$, $Q_{m}$ and $Q_\infty$ 
on the same probability space, in such a way that the equalities
$$B_{\alpha n^{1/4}}(Q_{n}) = B_{\alpha n^{1/4}}(Q_{m}) =  B_{\alpha n^{1/4}}(Q_{\infty})$$
  hold with probability at least $1-\varepsilon$.
\end{proposition}

\proof Let $\varepsilon>0$. From the first assertion of Proposition \ref{comptree}, we can find $\delta >0$ and an integer $n_0\geq 0$ such that
the following holds. If $n\geq n_0$ and $m\in\{n+1,n+2,\ldots\}$ are fixed, we can construct on the same probability space
a uniformly distributed pointed labeled tree with $n$ edges $(T_{n},\xi_n,\ell_n)$, 
a uniformly distributed pointed labeled tree with $m$ edges $(T_{m},\xi_m,\ell_m)$, and a uniform  infinite labeled tree $\Theta_{\infty} = (T_{\infty},\ell_{\infty})$,
in such a way that the event
\begin{eqnarray*}
&&E_{m,n} =\left\{\mathscr{P}(T_{\infty}, \lfloor \delta n^{1/2} \rfloor) = 
\mathscr{P}((T_{m}, \xi_{m}), \lfloor \delta n^{1/2} \rfloor)= \mathscr{P}((T_{n}, \xi_{n}), \lfloor \delta n^{1/2} \rfloor)\ ,\right.\\
&&\qquad\qquad\ \left.\ell_{\infty}\Big|_{\mathscr{P}(T_{\infty}, \lfloor \delta n^{1/2} \rfloor)} 
= \ell_{m}\Big|_{\mathscr{P}((T_{m}, \xi_{m}), \lfloor \delta n^{1/2} \rfloor)}= \ell_{n}\Big|_{\mathscr{P}((T_{n}, \xi_{n}), \lfloor \delta n^{1/2} \rfloor)} \right\}
\end{eqnarray*}
holds with probability at least $1-\varepsilon/2$. 
 Set
  \begin{eqnarray*} M_{n} &=& - \min_{0 \leq i \leq \lfloor \delta n^{1/2} \rfloor}  \ell_\infty( \mathrm{S}_{T_\infty}(i)). \end{eqnarray*}
Since $M_{n}$ is distributed as the maximal value of a random walk started at $0$ with increments uniformly distributed over $\{-1,0,+1\}$ and stopped after $\lfloor \delta  n^{1/2} \rfloor$ steps, Donsker's invariance principle shows that there exists 
a constant $\delta'>0$, which does not depend on $n$, such that the event $F_{n}:=\{M_{n} \geq \delta' n^{1/4}\}$ holds with probability at least $1-\varepsilon/2$.  

Let $\eta$ be a Bernoulli variable of parameter $1/2$ independent of the triplet $(\Theta_{n},\Theta_m,\Theta_\infty)$. 
We  set $Q_{n} = \Phi((T_n,\ell_n), \eta)$, $Q_{m}=\Phi((T_m,\ell_m),\eta)$ and $Q_{\infty} = \Phi'(\Theta_{\infty}, \eta)$. By Lemma \ref{geometric}, the equalities 
$$B_{\frac{1}{2}(\lfloor \delta' n^{1/4}\rfloor-3)}(Q_n)=B_{\frac{1}{2}(\lfloor \delta' n^{1/4}\rfloor-3)}(Q_\infty)=B_{\frac{1}{2}(\lfloor \delta' n^{1/4}\rfloor-3)}(Q_m)$$
hold on the event $E_{m,n} \cap F_{n}$ whose probability is at least $1- \varepsilon$. 
We just have to take $\alpha=\delta'/4$ to complete the proof. \endproof

\subsection{Proof of Theorem \ref{scaling-limit}}
 Let $(k_n)_{n\geq 1}$ be a sequence of non-negative real numbers converging to $\infty$ such that $k_n=o(n^{1/4})$ as $n \to \infty$.  We will prove simultaneously that, for every $r>0$, 
 \begin{eqnarray} B_r(k_n^{-1} \cdot \mathbf{Q}_n) & \xrightarrow[n\to\infty]{(d)} & B_r(\boldsymbol{\mathcal{P}}),\label{discreteconv1}\\
 B_r(k_n^{-1} \cdot \mathbf{Q}_\infty) & \xrightarrow[n \to \infty]{(d)} & B_r(\boldsymbol{\mathcal{P}}),
 \label{discreteconv2} \end{eqnarray}
 where the convergence holds in distribution in $\mathbb{K}$. Both assertions of Theorem \ref{scaling-limit} follow from these convergences
 (see also the comments following the statement of the theorem).
 
 We take $r=1$ to simplify notation. We first notice that when proving  \eqref{discreteconv1}, we may replace $B_1(k_n^{-1} \cdot \mathbf{Q}_n)$
 by $k_n^{-1}\cdot B_{k_n}(Q_n)$, simply because the Gromov-Hausdorff distance between $B_1(k_n^{-1} \cdot \mathbf{Q}_n)=k_n^{-1}\cdot B_{k_n}(
 \mathbf{Q}_n)$
 and $k_n^{-1}\cdot B_{k_n}(Q_n)$ is bounded above by $1/k_n$ -- see the beginning of this section. Similarly, when proving \eqref{discreteconv2}, we 
 may replace $B_1(k_n^{-1} \cdot \mathbf{Q}_\infty)$
 by $k_n^{-1}\cdot B_{k_n}(Q_\infty)$.
 
 Let $\varepsilon>0$. By Proposition \ref{couplingBMBP}, we may find a constant $\lambda_0>0$ such that, for every 
 $\lambda\geq\lambda_0$, 
 we can construct the Brownian plane $\boldsymbol{ \mathcal{P}}$ and the Brownian map $\mathbf{m}_\infty$
 simultaneously on the same probability space, in such a way that
 the equality
 \begin{eqnarray} \label{egalite2} B_1(\lambda\cdot \mathbf{m}_\infty)  &=& B_1( \boldsymbol{ \mathcal{P}})  \end{eqnarray} 
 holds with probability at least $1-\varepsilon$. 

 Then choose $\alpha>0$ such that the conclusion of Proposition \ref{compmap} holds. We may assume that $\alpha<(2\lambda_0)^{-1}$. Without loss of generality we may also assume that $k_n \leq \alpha \lfloor n^{1/4}\rfloor$ for every $n$. We set $m_n = \lceil \alpha^{-1} k_n\rceil ^4$. Notice
 that $m_n\leq n$ and $k_n\leq \alpha m_n^{1/4}$. 
 Using Proposition \ref{compmap} (with $n$ replaced by $m_n$ and $m$ replaced by $n$) and the notation of this proposition, we can 
 for every sufficiently large $n$ construct $Q_n$, $Q_{m_n}$ and $Q_\infty$ simultaneously on the same
 probability space, in such a way that the equalities
 \begin{eqnarray} \label{egalite}  {B}_{k_n}(Q_n)= {B}_{k_n}(Q_\infty) = {B}_{k_n}(Q_{m_n}),  \end{eqnarray} 
hold with probability at least $1-\varepsilon$. From the convergence of uniform quadrangulations towards the Brownian map \eqref{cvbm}, we have
  \begin{eqnarray*}   (V({Q}_{m_n}), (\alpha^{-1}k_n)^{-1}\op{d}_{\rm gr},\rho_{(m_n)}) & \xrightarrow[n\to\infty]{(d)} & \left(\frac{8}{9}\right)^{1/4} \cdot \mathbf{m}_\infty.\end{eqnarray*}
  where $\rho_{(m_n)}$ is the root vertex of $Q_{m_n}$ as in \eqref{cvbm}. This implies  \begin{eqnarray}
   k_n^{-1}\cdot B_{k_n}(Q_{m_n})& \xrightarrow[n\to\infty]{(d)} & B_1\left(\lambda \cdot \mathbf{m}_\infty\right). \label{eq:1}
  \end{eqnarray}
where $\lambda:=\left(\frac{8}{9}\right)^{1/4}\alpha^{-1}$. Notice that $\lambda\geq \lambda_0$ since $\alpha<(2\lambda_0)^{-1}$. Hence, as already noted, we can construct 
the Brownian plane $\boldsymbol{ \mathcal{P}}$ and the Brownian map $\mathbf{m}_\infty$
 simultaneously on the same probability space in such a way that \eqref{egalite2} holds with probability at least $1-\varepsilon$.
Thus, for any bounded continuous function $F: \mathbb{K}\to \R$, we have 
  \begin{eqnarray*}E\left[\left|F\big( k_n^{-1}\cdot B_{k_n}(Q_n)\big)-F\big(B_{1}( \boldsymbol{ \mathcal{P}})\big)\right|\right]& \leq & 
  E\left[\left|F\big( k_n^{-1}\cdot B_{k_n}(Q_n)\big)-F\big( k_n^{-1}\cdot B_{k_n}(Q_{m_n})\big)\right|\right] \\ &+&    
  E\left[\left| F\big( k_n^{-1}\cdot B_{k_n}(Q_{m_n})\big) - 
  F\big(B_1(\lambda \cdot \mathbf{m}_{\infty})\big)\right|\right] \\ &+&   
    E\left[\left|F\big(B_1(\lambda\cdot \mathbf{m}_{\infty})\big)-F\big(B_{1}( \boldsymbol{ \mathcal{P}})\big)\right|\right].  \end{eqnarray*}
By \eqref{egalite} and \eqref{egalite2},  the first and the third terms in the right-hand side of the last display are bounded by $2\varepsilon\, \sup|F|$, whereas  the convergence \eqref{eq:1} entails that the second term tends to $0$ as $ n \to \infty$. We 
have proved that  $ k_n^{-1}\cdot B_{k_n}(Q_n) \to B_{1}( \boldsymbol{ \mathcal{P}})$ in distribution in the Gromov-Hausdorff sense, and as noted at the beginning of the proof, this suffices to get
\eqref{discreteconv1}.  The proof of \eqref{discreteconv2} is similar: just replace $B_{k_n}(Q_n) $
by $B_{k_n}(Q_\infty)$ in the last display, still using \eqref{egalite}. This completes the proof of the theorem.\endproof

\subsection{An application to separating cycles} 

Proposition \ref{compmap} can be used to relate the large scale properties of the UIPQ to those of  random quadrangulations 
and of the Brownian map. The sphericity of the Brownian map was already applied in \cite{LGP} to prove the non-existence of small bottlenecks in large uniform quadrangulations. We present here a similar property of the UIPQ,  which partially answers a question raised by Krikun \cite{Kri}. \medskip 

Recall from subsection \ref{disquad} our notation $\op{Ball}_r(Q)$ for the ``combinatorial ball'' of radius $r$ associated with
a finite or infinite quadrangulation $Q$. Notice that $\op{Ball}_r(Q)$ does not determine the metric ball 
$B_r(Q)$ (and the converse does not hold either) because the knowledge of $\op{Ball}_r(Q)$ does not determine 
distances between vertices of $B_r(Q)$. However, a minor modification in the proof of Lemma \ref{geometric}
shows that the assumptions of this lemma also imply that $\op{Ball}_{r'}(Q)=\op{Ball}_{r'}(Q')$ -- just notice that 
all edges of the submap $\op{Ball}_{r'}(Q)$ will start and end at a corner of the pruned tree  $\mathscr{P}((\tau,\xi),h)$. 
Hence the statement of Proposition \ref{compmap} also remains valid if we replace the metric balls 
$B_{\alpha n^{1/4}}(Q_n)$, $B_{\alpha n^{1/4}}(Q_m)$ and $B_{\alpha n^{1/4}}(Q_\infty)$ by the combinatorial balls
of the same radius.

Let $Q$ be a finite or infinite rooted quadrangulation. Let $p\geq 2$, let $e_1,\ldots,e_p$ be \emph{oriented} edges of  $Q$
and let $x_1,\ldots,x_p$ be the respective origins of $e_1,\ldots,e_p$. Also set $x_{p+1}=x_1$. We say that $(e_1,\ldots,e_p)$ is an injective cycle of $Q$
of length $p$
if $x_1,\ldots,x_p$ are distinct and if $x_{i+1}$ is the target of $e_i$, for every $1\leq i\leq p$. When $p=2$,
we also exclude the case when $e_2$ is the same edge as $e_1$ with reverse orientation. If 
$\mathcal{C}=(e_1,\ldots,e_p)$ is an injective cycle of $Q$, the concatenation of $e_1,\ldots, e_p$ gives a closed loop
on the sphere, whose complement has exactly two connected components by Jordan's theorem. Suppose that 
$Q$ is infinite. We say that 
the cycle $\mathcal{C}$ separates the origin from infinity if the root vertex of $Q$ is not a vertex of $\mathcal{C}$ and if 
the connected component containing the root vertex contains finitely many vertices of $Q$.
In what follows, we say cycle instead of injective cycle. 

Consider the special case when $Q$ is the UIPQ $Q_\infty$. For every integer $n\geq 2$, Krikun \cite[Section
   3.5]{Kri} constructs a cycle $\mathcal{C}_{n}$ separating the origin
 from infinity in $Q_\infty$, which is contained in $B_{2n}(Q_\infty)\backslash B_n(Q_\infty)$, and such that its expected 
 length grows linearly in $n$ when $n\to\infty$. Krikun \cite[Conjecture 1]{Kri} also conjectures that the minimal length of  a
 cycle contained in the complement of $B_n(Q_\infty)$ and separating the origin from infinity must grow linearly 
 in $n$, a.s. The following corollary is a first step towards this conjecture.

 \begin{corollary}
 \label{Krikunconjec}
  Let $\kappa>1$ be an integer and let $f : \mathbb{N} \to \mathbb{R}_{+}$ be a function such that $f(n)=o(n)$ as $n\to \infty$. The probability that there exists an injective cycle of $Q_\infty$ with length less than $f(n)$, which separates the origin  from infinity 
 and whose vertices belong to $B_{\kappa n}(Q_\infty)\backslash B_n(Q_\infty)$, tends to $0$ as $n \to\infty$.
 \end{corollary}
 
 \proof Fix $\varepsilon >0$, and choose $\alpha>0$ small enough so that the conclusion of Proposition \ref{compmap} 
 (or rather of the variant of this proposition for combinatorial balls, as explained above) holds
 for every sufficiently large $n$.
Set   \begin{eqnarray*}m &=& \lceil\left( \frac{2\kappa n}{\alpha}\right)^4 \rceil \end{eqnarray*} 
 and notice that $\alpha m^{1/4} \geq 2\kappa n$. 
If $n$ is sufficiently large, we can then construct the UIPQ $Q_{\infty}$ and a uniform rooted quadrangulation  $Q_{m}$ with $m$ faces, in such a way that
 the equality
\begin{equation}
\label{equalityballs}
\op{Ball}_{2\kappa n}(Q_\infty)= \op{Ball}_{2\kappa n}(Q_m)
\end{equation}
holds  with probability at least $1-\varepsilon$. If the event considered in the proposition holds, there exists
a cycle $\mathcal{C}$ of $Q_\infty$
with length less than $f(n)$ and whose vertices belong to $B_{\kappa n}(Q_\infty)\backslash B_n(Q_\infty)$, and a vertex $v$ of $Q_\infty$ at graph distance $2\kappa m$
from the root vertex, such that $v$ and the root vertex of $Q_\infty$ belong to distinct components of the complement of
the cycle (choose $v$ such that there exists a path from $v$ to infinity that stays outside $B_{2\kappa n-1}(Q_\infty)$). 
Assuming that the identity (\ref{equalityballs}) holds, we see that on the event
of the proposition there exists a cycle $\mathcal{C}'$ of $Q_m$, with length less than $f(n)$ and  whose vertices belong to $B_{\kappa n}(Q_m)\backslash B_n(Q_m)$,
and a vertex $v'$ at graph distance $2\kappa m$
from the root vertex of $Q_m$, such that $v'$ and the root vertex of $Q_m$ belong to distinct components of the complement of
the cycle $\mathcal{C}'$. Now note that the set of all vertices of $Q_m$ lying in the connected component containing $v'$ must have diameter 
at least $\kappa n$ for the graph distance, 
whereas the set of all vertices lying in the connected component of the root vertex must have diameter at least $n$. By Corollary 1.2 in \cite{LGP}, we get
that the probability that both (\ref{equalityballs}) and the event of the proposition hold tends to $0$ as 
$n\to\infty$. This completes the proof. 
\endproof

\section{Properties of the Brownian plane}
\label{sec:properties}

\subsection{Topology and Hausdorff dimension of the Brownian plane}

The following proposition is analogous to a result of \cite{IM}. It gives a more explicit form 
to the definition of the equivalence relation $\approx$ in Section 1.

\begin{proposition}
\label{identrel}
Almost surely, for every $a,b\in\t_\infty$, the property $a\approx b$ holds 
if and only if $D^\circ_\infty(a,b)=0$.
\end{proposition}

\proof By definition, $a\approx b$ if and only if $D_\infty(a,b)=0$. Since 
$D_\infty(a,b)\leq D^\circ_\infty(a,b)$, it is obvious that 
the property $D^\circ_\infty(a,b)=0$ implies $a\approx b$. We need 
to prove the converse. We fix $r>0$ and $\delta\in(0,1)$. It will be enough
to prove that, with probability at least $1-\delta$, for every 
$a,b\in\t_\infty$ such that $D_\infty(\rho_\infty,a)\leq r$ and 
$D_\infty(\rho_\infty,b)\leq r$, the property $D_\infty(a,b)=0$ 
implies $D^\circ_\infty(a,b)=0$.

To this end, we rely on the coupling argument of Proposition \ref{couplingBMBP}, and we use 
the notation of the proof of this proposition. In particular, we fix 
$\lambda \geq \lambda_0$, and we argue on the event $\mathcal{F}_\lambda$ (recall that the probability of this event
is bounded below by $1-\delta$). We also use the notation $T=\lambda^4$. Let $a,b\in\t_\infty$ such that $D_\infty(\rho_\infty,a)\leq r$ and 
$D_\infty(\rho_\infty,b)\leq r$, and assume that $D_\infty(a,b)=0$. Write $a=p_\infty(s')$ and $b=p_\infty(t')$ for some 
$s',t'\in\R$. By Lemma \ref{keyasymp} (ii), we must have 
$s',t'\in (-\eta_\infty(A),\gamma_\infty(A))$. Set $s=s'$ if $s'\geq 0$ and
$s=T+s'$ if $s'<0$ and define $t$ similarly from $t'$. Note that
$s,t\in [0,\gamma_\lambda(A))\cup (\eta_\lambda(A),T]$. By Lemma \ref{corokey},
$$D^*_{\lambda}(p_{(\lambda)}(s),p_{(\lambda)}(t))=D_\infty(p_\infty(s'),p_\infty(t'))=0.$$
By \cite[ Theorem 3.4]{IM}, the condition $D^*_{\lambda}(p_{(\lambda)}(s),p_{(\lambda)}(t))=0$
implies $D^\circ_{\lambda}(p_{(\lambda)}(s),p_{(\lambda)}(t))=0$.

Consider first the case when $p_{(\lambda)}(s)=p_{(\lambda)}(t)$, meaning that
$$\be^\lambda_s=\be^\lambda_t = \min_{s\wedge t\leq r\leq s\vee t} \be^\lambda_r.$$
Since $s,t\in [0,\gamma_\lambda(A))\cup (\eta_\lambda(A),T]\subset[0,\alpha]\cup[T-\alpha,T]$, and
we know that
$$\min_{\alpha\leq r\leq T-\alpha} \be^\lambda_r=\inf_{r\geq \alpha} R_r \wedge \inf_{r\geq \alpha} R'_r > A^4,$$
it easily follows that $X_{s'}=X_{t'}=m_X(s',t')$, and consequently $a=p_\infty(s')=p_\infty(t')=b$, so that
obviously $D^\circ_\infty(a,b)=0$.

Consider then the case when $p_{(\lambda)}(s)\not =p_{(\lambda)}(t)$. Then, the fact
that $D^\circ_{\lambda}(p_{(\lambda)}(s),p_{(\lambda)}(t))=0$ implies that $p_{(\lambda)}(s)$ and
$p_{(\lambda)}(t)$ are both leaves of $\t_{(\lambda)}$ (see Lemma 3.2 in \cite{LGP}). It follows that
$$D^\circ_{\lambda}(s,t)= D^\circ_{\lambda}(p_{(\lambda)}(s),p_{(\lambda)}(t))=0.$$
Hence we  have
\begin{equation}
\label{identech1}
\min_{r\in[s\wedge t,s\vee t]} Z^{\lambda}_r = Z^{\lambda}_s = Z^{\lambda}_t
\end{equation}
or 
\begin{equation}
\label{identech2}
\min_{r\in[0,s\wedge t]\cup [s\vee t,T]} Z^{\lambda}_r = Z^{\lambda}_s = Z^{\lambda}_t.
\end{equation}
Suppose first that $s,t\in [0,\gamma_\lambda(A))$. Then since
$$\min_{r\in[\gamma_\lambda(A),\eta_\lambda(A)]} Z^\lambda_r < -10 r$$
and $Z^{\lambda}_s=Z_{s'}\geq -r$ (because $D_\infty(\rho_\infty,a)\leq r$), it is
obvious that  (\ref{identech2}) cannot hold, so that (\ref{identech1}) holds.
Similarly, if $s,t\in (\eta_\lambda(A),T]$, we get that (\ref{identech1}) holds.
Finally, if $s\in[0,\gamma_\lambda(A))$ and $t\in(\eta_\lambda(A),T]$ (or conversely), 
we obtain that (\ref{identech2}) holds.

In all three cases, we can now verify that
$$\min_{r\in [s'\wedge t',s'\vee t']} Z_r = Z_{s'} = Z_{t'}$$
so that $D^\circ_\infty(s',t')=0$, and $D^\circ_\infty(a,b)=0$. This completes the proof. \endproof

\begin{proposition}
\label{Hausdorffplane}
The Hausdorff dimension of $\pp$ is almost surely equal to 
$4$.
\end{proposition}

This statement immediately follows from Proposition \ref{couplingBMBP}, together with the fact that the Hausdorff dimension of the
Brownian map, or of any nontrivial ball in the Brownian map, is equal to $4$,
see Theorem 6.1 in \cite{IM}. We leave the details to the reader.

\begin{proposition}
\label{homeoplane}
Almost surely, the space $\pp$ is homeomorphic to the plane.
\end{proposition}

\proof We adapt the arguments of \cite{LGP} to our setting. We write $\widehat\pp$
for the Alexandroff compactification of $\pp$. We will prove that 
$\widehat\pp$ is almost surely homeomorphic to the sphere $\SS^2$.
Recall the definition of the equivalence relation $\sim_X$
in the first section, and with a slight abuse of notation, define 
the relation $\approx$ on $\R$ by setting $s\approx t$ if
and only if
$$Z_s=Z_t = \min_{r\in[s\wedge t,s\vee t]} Z_r$$
or equivalently if $D^\circ_\infty(s,t)=0$. Note that $\approx$ is also
an equivalence relation on $\R$. Furthermore, it follows from Lemma
3.2 in \cite{LGP} (and the coupling argument of the proof of Proposition \ref{couplingBMBP})
that, almost surely for every $s,t,u\in\R$, the properties
$s\sim_X t$ and $s\approx u$ may hold simultaneously only 
if $s=t$ or $s=u$. It follows that, outside a set of probability zero
which we discard from now on, we can define another equivalence
relation $\simeq$ on $\R$ by setting $s\simeq t$ if and only if
$s\sim_X t$ or $s\approx t$. 

Write $\widehat\R=\R\cup\{\infty\}$ for the Alexandroff compactification of $\R$,
and extend both equivalence relations $\sim_X$ and $\simeq$ to
$\widehat\R$ by declaring that the equivalence class of $\infty$
is a singleton. Equip $\widehat\R/\!\sim_X$ with the quotient topology
(clearly the resulting space is the Alexandroff compactification of $\t_\infty$).
Let $\Pi: \R/\!\sim_X=\t_\infty\longrightarrow \pp$ be the canonical projection, and extend
it to a projection from
$\widehat \R/\!\sim_X$ onto  $\widehat \pp$ by
mapping the equivalence class of $\infty$ to the point at infinity in $\widehat\pp$. Then $\Pi$ 
is continuous. The only point to check is the continuity at $\infty$, which follows from
the easy fact that $\Pi^{-1}(B_r(\boldsymbol{ \mathcal{P}}))$ is a compact subset of $\t_\infty$, for every $r\geq 0$. 
The projection $\Pi$ then factorizes through a bijection from $\widehat\R/\!\simeq$
onto $\widehat\pp$, which is also continuous hence provides a homeomorphism
from $\widehat\R/\!\simeq$ onto $\widehat \pp$. 

We then use Proposition 2.4 in \cite{LGP} to see that the quotient space
$\widehat\R/\!\simeq$ is a.s. homeomorphic to the sphere $\SS^2$. Our setting is slightly
different since \cite{LGP} deals with a quotient space of the circle $\SS^1$. This
makes no real difference since $\widehat\R$ is homeomorphic to $\SS^1$. 
Another difference is the fact that the random functions $X$ and $Z$
used to define the equivalence relations $\sim_X$ and $\approx$
are not defined at the point $\infty$, whereas \cite{LGP} considers 
functions that are (continuous and) defined everywhere on the circle.
Nonetheless one can easily check that the arguments of the proof of 
Proposition 2.4 in \cite{LGP} go through without change, provided we
can verify that the local minima of $X$, respectively of $Z$, are distinct.
In the case of $X$, this is a standard fact that follows from the connections
between linear Brownian motion and the three-dimensional Bessel process.
The case of $Z$ is treated by arguments similar to the proof of Lemma 3.1
in \cite{LGP}.

We can now conclude that $\widehat \pp$ is (almost surely) homeomorphic 
to $\SS^2$, and the statement of the proposition follows. \endproof

\subsection{Geodesic rays in the Brownian plane}

The results of this section are analogous to the discrete results proved in
\cite[Section 3.2]{CMM} for the uniform infinite planar quadrangulation, and are also
closely related to the study of geodesics in the Brownian map \cite{AM}.

Let $x\in \pp$. A geodesic ray from $x$ is an infinite continuous path $\omega:[0,\infty)\longrightarrow \pp$
such that $\omega(0)=x$ and $D_\infty(\omega(s),\omega(t)) = |t-s|$ for every $s,t\geq 0$. 

Recall our notation $\bp_\infty$ for the canonical projection from $\R$ onto $\pp$.
Fix $x\in\pp$ and let $t\in\R$ such that $\bp_\infty(t)=x$. For every $r\geq 0$, set
$$\gamma_t(r)=\inf\{s\geq t : Z_s=Z_t-r\}.$$
It is clear that $\gamma_t(r)<\infty$ for every $r\geq 0$, a.s. Then set
$\omega_t(r)=\bp_\infty(\gamma_t(r))$,
for every $r\geq 0$. 

\begin{proposition}
\label{simpleray}
For every $x\in\pp$ and $t\in\R$ such that $\bp_\infty(t)=x$, $\omega_t$ is a geodesic
ray from $x$. Such a geodesic ray will be called a simple geodesic ray.
\end{proposition}

\proof
It is obvious that $\gamma_t(0)=t$ and thus $\omega_t(0)=x$.
Then, by definition $Z_{\gamma_t(r)}= Z_t-r$, for every $r\geq 0$. It follows that,
for every $r,s\geq 0$,
$$D_\infty(\omega_t(r),\omega_t(s)) \geq |Z_{\gamma_t(r)} - Z_{\gamma_t(s)}| =|r-s|.$$
On the other hand, it is also immediate from the definition of $\gamma_t(r)$, that,
for every $r,s\geq 0$,
$$D^\circ_\infty(\gamma_t(r),\gamma_t(s))=|r-s|.$$
The equality $D_\infty(\omega_t(r),\omega_t(s)) =|r-s|$ now follows. \endproof

\begin{proposition}
\label{unirayorigin}
Almost surely, $\omega_0$ is the unique geodesic ray from $\rho_\infty$.
\end{proposition}

\proof By Proposition \ref{couplingBMBP}, 
for every $\delta>0$, we can find $\lambda>0$ large enough so that, with probability at
least $1-\delta$, the ball $B_1(\boldsymbol{ \mathcal{P}})$
is isometric to the ball $B_1(\lambda\cdot \bm_\infty)$. By Corollary 7.7 in \cite{AM} (and the invariance of the Brownian map
under re-rooting, see Theorem 8.1 in \cite{AM}), we can find a (random) $\ve>0$
such that all geodesic paths from $\rho$ to a point outside $B_1(\lambda\cdot \bm_\infty)$
coincide over the interval $[0,\ve]$. Now, let ${\mathcal R}$ be the set of all
geodesic rays from $\rho_\infty$, and set
$$\tau = \inf\{t\geq 0: \exists \omega\in {\mathcal R} : \omega(t)\not =\omega_0(t)\}.$$
By the previous considerations, $\tau >0$ a.s. However, the scaling invariance of 
$\boldsymbol{ \mathcal{P}}$ guarantees that $\tau$ has the same distribution as $\lambda\tau$, for
every $\lambda >0$. It follows that $\tau=\infty$ a.s., which completes the proof. \endproof

Our goal is to prove that all geodesic rays in $\pp$ are simple geodesic rays. We will
rely on the preceding proposition and on the following invariance property of
the Brownian plane under re-rooting.

\begin{proposition}
\label{re-root}
Let $t\in \R$. The pointed space $(\pp,D_\infty,\bp_\infty(t))$ has the same
distribution as $\boldsymbol{ \mathcal{P}}=(\pp,D_\infty,\rho_\infty)$.
\end{proposition}

\proof Suppose that $t>0$ for the definiteness. Define two processes
$\widetilde R$ and $\widetilde R'$ by setting, for every $s\geq 0$,
$$\widetilde R_s = R_t + R_{t+s} - 2\min_{t\leq r\leq t+s} R_r$$
and
$$\widetilde R'_s=\left\{
\begin{array}{ll}
R_t+R_{t-s} - 2\min_{t-s\leq r\leq t} R_r\quad&\hbox{if }s\leq t,\\
R'_{s-t}+R_t -2\min_{r\in(-\infty,t-s]\cup [t,\infty)} X_r&\hbox{if }s>t.
\end{array}
\right.$$
With the notation of Section 1, we have $\widetilde R_s=d_X(t,t+s)$ and $\widetilde R'_s=d_X(t,t-s)$ for every $s\geq 0$.

We also set $\widetilde Z_s=Z_{t+s}-Z_t$ for every $s\in\R$. Then, we claim that the two triplets $(\widetilde R,\widetilde R',\widetilde Z)$
and $(R,R',Z)$
have the same distribution. In the case when $\t_\infty$ is replaced by the CRT, a similar ``re-rooting invariance'' identity
can be found as Corollary 4.9 in Marckert and Mokkadem \cite{MaMo}, and our claim then follows from a suitable passage to
the limit: Just apply the Marckert-Mokkadem result to the (scaled) CRT coded by a Brownian excursion of
duration $T$ and let $T$ tend to $\infty$ in the spirit of Section 3 above. Alternatively, the reader who is
familiar with properties of the three-dimensional Bessel process will be able to verify that the pairs 
$(R,R')$ and $(\widetilde R,\widetilde R')$ have the same distribution, from which our claim also follows
in a straightforward manner. Notice that, in the same way as we defined $\t_\infty$ in Section 1, we 
can associate a random tree $\widetilde \t_\infty$ with the pair $(\widetilde R,\widetilde R')$, and 
$\widetilde \t_\infty$ is easily identified to $\t_\infty$ ``re-rooted'' at $p_\infty(t)$.

To complete the proof, just note that the pointed space 
$(\pp,D_\infty,\bp_\infty(t))$ can be obtained from the triplet $(\widetilde R,\widetilde R',\widetilde Z)$ by the same
construction that we used to obtain $(\pp,D_\infty,\rho_\infty)$ from
$(R,R',Z)$. \endproof

\rem By combining the last proposition with the preceding one, we obtain
that almost surely for every rational $t$ the simple geodesic ray $\omega_t$
is the unique geodesic ray from $\bp_\infty(t)$. This will be useful in
the proof of Theorem \ref{uniray} below.

\medskip
The next proposition shows that the quantities $Z_a$, $a\in\t_\infty$ can be
interpreted as measuring the relative distances from the point at infinity
in the Brownian plane. Recall that $\Pi$ stands for the canonical projection from
$\t_\infty$ onto $\pp$. The following proposition should be compared to \cite[Theorem 1]{CMM}.

\begin{proposition}
\label{Zdistance}
Almost surely, for every $a,b\in\t_\infty$,
$$Z_a - Z_b \ =\ \lim_{x\to\infty} ( D_\infty(\Pi(a),x)- D_\infty(\Pi(b),x))$$
where the limit holds when $x$ tends to the point at infinity in the Alexandroff
compactification of $\pp$. Consequently, if $\omega$
is any geodesic ray in $\pp$, we have
$$Z_{\omega(r)}=Z_{\omega(0)}-r$$
for every $r\geq 0$.
\end{proposition}

\proof We fix $s,t\in\R$ and take $a=p_\infty(s)$ and $b=p_\infty(t)$. To
simplify notation, we set
$$m_Z(s,t)=\min_{r\in[s\wedge t,s\vee t]} Z_r.$$
From our construction, it is clear that the simple geodesic rays 
$\omega_s$ and $\omega_t$ coalesce in finite time. More
precisely, we have for every $r\geq 0$,
\begin{equation}
\label{Zdistech}
\omega_s(Z_s-m_Z(s,t) + r) = \omega_t(Z_t-m_Z(s,t)+r).
\end{equation}
For every integer $n\geq 1$, set
$$S_n(a)=\{x\in \pp : D_\infty(\Pi(a),x)=n\}.$$
Let ${\rm Geo}_n(a)$ be the set of all geodesic paths 
from $\Pi(a)$ to a point of $S_n(a)$, and 
$$\eta_n(a)=\inf\{r>0: \exists \omega,\omega'\in {\rm Geo}_n(a) : \omega(r)\not =\omega'(r)\}.$$
By combining the invariance of the Brownian plane under re-rooting (Proposition
\ref{re-root}) with the argument already used in the proof of Proposition \ref{unirayorigin},
we get that $\eta_n(a)>0$ a.s. On the other hand, the sequence $(\eta_n(a))_{n\geq 1}$
is clearly increasing and if $\eta_\infty(a)$ denotes its limit, the scale invariance 
of the Brownian plane implies that $\lambda\eta_\infty(a)$ has the same distribution as 
$\eta_\infty(a)$, for every $\lambda >0$. It follows that $\eta_\infty(a)=\infty$ a.s.

Consequently, we can choose $n$ sufficiently large so that $\eta_n(a)> Z_s-m_Z(s,t)$.
Then, let $x\in \pp$ such that $D_\infty(\Pi(a),x)\geq n$, and let $\bar\omega$ be any geodesic path
from $\Pi(a)$ to $x$ (such a geodesic exists because the 
Brownian plane is a boundedly compact length space, 
see subsection 2.1). By the definition of $\eta_n(a)$, we have $\bar\omega(r)=\omega_s(r)$
for every $r\in[0,\eta_n(a)]$. Recalling (\ref{Zdistech}), we can construct a continuous 
path $\omega'$ from $\Pi(b)$ to $x$ by concatenating the paths
$$(\omega_t(r),0\leq r \leq Z_t-m_Z(s,t))$$
and
$$(\bar\omega(Z_s-m_Z(s,t)+r),0\leq r\leq D_\infty(\Pi(a),x) - (Z_s-m_Z(s,t))).$$
The length of the path $\omega'$ is
$$Z_t-m_Z(s,t) + (D_\infty(\Pi(a),x) - (Z_s-m_Z(s,t))) = D_\infty(\Pi(a),x) + (Z_t-Z_s)$$
and so we have obtained the bound
$$D_\infty(\Pi(b),x) \leq D_\infty(\Pi(a),x) + (Z_t-Z_s)$$
which holds for any $x\in\pp$ such that $D_\infty(\Pi(a),x)\geq n$.
It follows that
$$\liminf_{x\to\infty} ( D_\infty(\Pi(a),x)- D_\infty(\Pi(b),x)) \geq Z_s -Z_t =Z_a -Z_b.$$
We can now interchange the roles of $a$ and $b$ and we get that the convergence
of the proposition holds for this particular choice of $a$ and $b$, almost surely.

Finally, the convergence of the proposition holds  outside a set of probability zero
for all $a,b$ of the form $a=p_\infty(s)$ and $b=p_\infty(t)$ with $s,t\in\Q$. A simple
density argument now completes the proof of the first assertion. 

The last assertion of the proposition is then immediate: If $\omega$ is any geodesic ray, we have for every $0\leq r\leq r'$,
$$D_\infty(\omega(r),\omega(r'))-D_\infty(\omega(0),\omega(r'))=-r$$
and $\omega(r')\longrightarrow \infty$ as $r'\to\infty$, so that we just need to apply the first assertion. 
\endproof

\rem The preceding proof shows that, if $a=p_\infty(s)$ and $b=p_\infty(t)$ for some fixed
$s,t\in\R$, we have indeed $D_\infty(\Pi(a),x)- D_\infty(\Pi(b),x)=Z_a - Z_b$
for all $x\in\pp$ such that $D(\rho_\infty,x)$ is sufficiently large, almost surely.

\begin{theorem}
\label{uniray}
All geodesic rays in $\pp$ are simple geodesic rays.
\end{theorem}

In particular, this implies that any pair of geodesic rays coalesces in finite time.

\proof Let $x\in\pp$, and let 
$\omega$ be a geodesic ray from $x$. We first assume that $x=\Pi(a)$
where $a$ is a leaf of $\t_\infty$ (the case when $a$ is not a leaf 
will be treated at the end of the proof). Then, there is a unique
$t\in\R$ such that $a=p_\infty(t)$. Fix $u>0$. We will prove that
$\omega(u)=\omega_t(u)$. Recall the notation
$$\gamma_t(u)=\inf\{s\geq t : Z_s=Z_t-u\},$$
and also set
$$\gamma'_t(u)=\sup\{s\leq t : Z_s=Z_t-u\}.$$
Note that $\bp_\infty(\gamma_t(u))= \bp_\infty(\gamma'_t(u))$, and
that $p_\infty(\gamma_t(u))$ and  $p_\infty(\gamma'_t(u))$ are both leaves
of $\t_\infty$ (the latter property is a consequence of the fact that, for $t_1,t_2,t_3\in \R$,
the properties $t_1\sim_X t_2$ and $t_1\approx t_3$ may hold simultaneously only if 
$t_1=t_2$ or $t_1=t_3$, as was mentioned in the proof of Proposition \ref{homeoplane}). Let $c_0$ be the unique vertex of $\t_\infty$ such that
$$\llbracket a,p_\infty(\gamma_t(u))\rrbracket \cap \llbracket a,p_\infty(\gamma'_t(u))\rrbracket=\llbracket a,c_0\rrbracket.$$
Fix a point $c_1\in\rrbracket a,c_0\llbracket$ such that the range of $\omega$ does not contain $\Pi(c_1)$. Such
a point exists because, for some values of $r\in\R$, the set $\rrbracket a,c_0\llbracket$ contains
uncountably many vertices $c$ such that $Z_c=r$, and, by the preceding proposition, the 
geodesic ray $\omega$ can visit $\Pi(c)$ for at most one such vertex $c$. We can also assume
that $c_1$ is not a branching point of $\t_\infty$, because the set of all branching points 
is countable. Then let $\t_1$ be the connected component of $\t_\infty\backslash\{c_1\}$
that contains $a$, and note that $\t_1$ is bounded (otherwise this would contradict the 
definition of $\gamma_t(u)$ and $\gamma'_t(u)$). 

Let $T=\inf\{s\geq 0: \omega(s)\notin \Pi(\t_1)\}$.
By a simple argument (see the beginning of Section 3 in \cite{AM}), there must
exist two vertices $b$ and $b'$ such that $b\in\t_1$, $b'\in\t_\infty\backslash\t_1$
and $\omega(T)=\Pi(b)=\Pi(b')$, so that in particular $b\approx b'$. Notice that
$b$ and $b'$ must be leaves of $\t_\infty$, and we can define $s_1,s'_1\in \R$
by the conditions $p_\infty(s_1)=b$, $p_\infty(s'_1)=b'$. Since $b\approx b'$,
we have $Z_s\geq Z_{s_1}=Z_{s'_1}$ for every $s\in[s_1\wedge s'_1,s_1\vee s'_1]$.
We can now pick any rational $s$ in $[s_1\wedge s'_1,s_1\vee s'_1]$, and we obtain from
our definitions that the range of the geodesic ray $\omega_s$
must contain $\bp_\infty(s_1\vee s'_1)= \omega(T)$. 
So there exists $u_1\geq 0$ such that $\omega_s(u_1)=\omega(T)$. Let $\widetilde \omega$
be the infinite path obtained from the concatenation of $(\omega_s(r),0\leq r\leq u_1)$
and $(\omega(r),r\geq T)$. Then, it
is easy to verify that $\widetilde \omega$ is a geodesic ray: If $r\in[0,u_1]$
and $r'\in[u_1,\infty)$, the bound $D_\infty(\widetilde \omega(r),\widetilde \omega(r'))\leq r'-r$
is clear from the triangle inequality, but the reverse bound is also easy by
writing $D_\infty(\widetilde \omega(r),\widetilde \omega(r'))\geq |Z_{\tilde\omega(r)}- Z_{\tilde\omega(r')}|$. 

Finally, since $\widetilde \omega$ is a geodesic ray starting from $\bp_\infty(s)$
with a rational value of $s$, we know that $\widetilde \omega$ must coincide
with the simple geodesic ray $\omega_s$. Since $\omega_s$ clearly visits 
$\omega_t(u)$ it follows that $\omega$ also visits $\omega_t(u)$, and finally
$\omega_t(u)=\omega(u)$.

It remains to consider the case when $a$ is not a leaf of $\t_\infty$. In that case, we can find
arbitrarily small values of $r>0$ such that $\omega(r)=\Pi(b)$, where $b$ is leaf of $\t_\infty$
(otherwise $\omega$ would have to visit an entire geodesic segment of $\t_\infty$, which
is absurd). By the first part of the proof, there are arbitrarily small values of $r>0$
such that $(\omega(r+u))_{u\geq 0}$ is a simple geodesic ray. The desired result
easily follows.
\endproof

\rem If $x=\Pi(a)$ and $a$ is not a leaf of $\t_\infty$, there are two (three if $b$
is a branching point) distinct geodesic rays starting from $x$. This should be
compared with Theorem 1.4 in \cite{AM}.

\end{document}